\documentclass[11pt]{article}
\usepackage{amsmath,epsfig,subfigure,amsthm,amsfonts,amssymb}
\usepackage{graphics, graphicx}
\usepackage[section]{placeins}
\usepackage{cases}
\usepackage[rightcaption]{sidecap}
\usepackage{multicol}
\usepackage[sort]{natbib}
\newtheorem{thm}{Theorem}[section]

\newtheorem{lem}[thm]{Lemma}

\newtheorem{rem}[thm]{Remark}
\newtheorem{defn}[thm]{Definition}

\newcommand{\thmref}[1]{Theorem~{\rm \ref{#1}}}
\newcommand{\lemref}[1]{Lemma~{\rm \ref{#1}}}

\makeatletter \@addtoreset{equation}{section}

\newcommand{\beq}[1]{\begin{equation} \label{#1}}
\newcommand{\eeq}{\end{equation}}
\newcommand{\bed}{\begin{displaymath}}
\newcommand{\eed}{\end{displaymath}}

\newcommand{\bedd}{\bed\begin{array}{l}}
\newcommand{\eedd}{\end{array}\eed}
\newcommand{\bea}{\bed\begin{array}{rl}}
\newcommand{\eea}{\end{array}\eed}

\newcommand{\disp}{\displaystyle}
\newcommand{\ad}{&\!\!\!\disp}
\newcommand{\aad}{&\disp}
\newcommand{\barray}{\begin{array}{ll}}
\newcommand{\earray}{\end{array}}

\newcommand{\diag}{{\rm{diag}}}
\def\para#1{\vskip 0.1\baselineskip\noindent{\bf #1}}
\def\qed{\strut\hfill $\Box$}

\def\({\left(}
\def\){\right)}

\def\qed{\hfill{$\Box$}}
\def\F{{\cal F}}

\newcommand{\U}{{\mathcal U}}
\newcommand{\Ll}{{\cal L}}
\def\cd{(\cdot)}
\def\M{{\cal M}}
\def\A{{\cal A}}
\def\B{{\cal B}}
\def\H{{\cal H}}
\def\bb{\mathfrak{B}}
\def\rr{{\mathbb R}}

\newcommand{\var}{\text{Var}}

\newcommand{\la}{\lambda}

\newcommand{\Ga}{\Gamma}
\newcommand{\dl}{\delta}
\newcommand{\ka}{\kappa}
\newcommand{\Dl}{\Delta}
\newcommand{\e}{\varepsilon}

\newcommand{\ar}{\rightarrow}
\newcommand{\lf}{\lfloor}
\newcommand{\rf}{\rfloor}

\newcommand{\nd}{\noindent}

\newcommand{\al}{\alpha}

\newcommand{\sg}{\sigma}
\newcommand{\pr}{{P}}
\newcommand{\ex}{{\mathbb E}}
\newcommand{\wdt}{\widetilde}

\newcommand{\wrt}{{with respect to }}
\newcommand{\lbar}{\overline}

\newcommand{\exh}{\mathbb{E}^{\phi,h,\dl}_{x,i,n}}
\newcommand{\exr}{\mathbb{E}^{\phi,h,\dl}_{x,i,r}}
\newcommand{\prh}{\pr^{\phi,h,\dl}_{x,i,n}}
\newcommand{\vah}{\var^{\phi,h,\dl}_{x,i,n}}
\def\para#1{\vskip 0.4\baselineskip\noindent{\bf #1}}

\begin{document}
\title{Approximation of A Class of Non-Zero-Sum Investment and Reinsurance Games for Regime-Switching Jump-Diffusion  Models}

\author{Trang Bui,\thanks{Department of Mathematics, Wayne State University, Detroit, MI 48202. Email: trang.bui@wayne.edu. The research was supported in part by
the
Air Force Office of Scientific Research under Grant FA9550-18-1-0268.}
\and Xiang Cheng,\thanks{Centre for Actuarial Studies,
Department of Economics, The University of Melbourne, VIC 3010, Australia,
cheng.x@unimelb.edu.au.}
\and Zhuo Jin,\thanks{Centre for Actuarial Studies,
Department of Economics, The University of Melbourne, VIC 3010, Australia,
zjin@unimelb.edu.au.}
\and George Yin\thanks{Department of Mathematics, Wayne State University, Detroit, MI 48202.
Email: gyin@math.wayne.edu. The research was supported in part by
the
Air Force Office of Scientific Research under Grant FA9550-18-1-0268.}}

\maketitle
\begin{abstract}
This work develops an approximation procedure for a class of non-zero-sum stochastic differential investment and reinsurance games between two insurance companies.
Both proportional reinsurance and excess-of loss reinsurance policies are considered. We develop numerical algorithms to obtain the Nash equilibrium by adopting the Markov chain approximation methodology and applying the dynamical programming principle for the nonlinear integro-differential Hamilton-Jacobi-Isaacs (HJI) equations. Furthermore, we establish the convergence of the approximation sequences and the approximation to the value functions. Numerical examples are presented to illustrate the applicability of the algorithms.

\vskip 0.25 in \nd{\bf Key Words.}  Stochastic control, non-zero-sum game, investment, reinsurance, regime switching.

\end{abstract}


\newpage
\section{Introduction}\label{sec:intro}

Insurers tend to accumulate relatively large amounts of cash or cash equivalents through the written
insurance portfolio. Investing the surplus in a financial market in order to pay future claims and
to avoid financial ruin becomes a natural choice. In terms of financial performance the investment
income allows significant pricing flexibility in underwriting to the insurers. The surplus is allowed
to be invested
in a  financial market in continuous time.

On the other hand, reinsurance has been considered as an effective risk management tool for insurance companies
to transfer their risk exposure to another commercial institution. The primary insurer pays the reinsurer a certain portion
of the premiums. In return, the reinsurer is obliged to share the risk of large claims with
the primary insurer. Proportional reinsurance and excess-of-loss reinsurance are two major types
of reinsurance strategies. With proportional reinsurance, the reinsurance company covers a fixed
percentage of losses. The fraction of risk shared by the reinsurance company is determined when the
reinsurance contract is sold. The other type of reinsurance policy is nonproportional reinsurance.
The most common nonproportional reinsurance policy is excess-of-loss reinsurance, where the
 primary insurance carrier (called cedent)
will pay all of the claims up to a predetermined amount (termed
retention level).

The optimal risk controls for an insurance corporation has been studied extensively since the classical collective risk model was introduced by \cite{Lundberg03}. The insurance companies can reduce or eliminate the risk of loss by involving in a reinsurance program and reinvesting in the stock market. Particularly, advanced stochastic control theory and dynamic programming
principle are widely used to design
the optimal reinsurance
and investment strategies of an insurance company in previous works.
\cite{Asmussen} used diffusion approximation to find an optimal policy balancing the risk and expected profits for a financial corporation with excess-of-loss reinsurance. \cite{Choulli01} investigated the case of excess-of-loss reinsurance for an insurance company facing a constant liability payment to maximize the expected present value of total future dividend pay-outs. \cite{HaldS04} studied the optimal proportional reinsurance policy to maximize the adjustment coefficient of the ruin probability with the variance premium principle; see also \cite{Browne95}, \cite{ChenLL10}, \cite{ZhangS12}, and \cite{ZhangMZ16}.

Recently, the extension of optimal investment and reinsurance problem (Nash equilibrium) in the context of stochastic differential games including zero-sum games and non-zero-sum games has been developed rapidly; see the existence of the Nash equilibrium of non-zero-sum stochastic differential  game with $N$ players over infinite time horizon in \cite{BF2000}. The existence of the Nash equilibrium of a non-zero-sum stochastic differential game between two insurance companies in \cite{Zeng}. \cite{LiuY13} studied a zero-sum stochastic differential reinsurance and investment game between two competing insurance companies under VaR constraints for the purpose of risk management. \cite{BensoussanSYY} investigated a class of non-zero sum stochastic differential game between two insurers by using the objectives of relative performance and obtained explicit solutions for optimal reinsurance and investment strategies. \cite{ChenS17} formulated a stochastic Stackelberg
differential reinsurance game between an insurer and a reinsurer, allowing them to consider the benefits of both parties in the reinsurance contracts.
Under the criteria of maximizing the expected utilities of the players' terminal surpluses, the Stackelberg equilibrium strategies and value functions are obtained by using the backward
stochastic differential equation (BSDE) approach; see also \cite{TaksarZ11}, \cite{PunW16}, and \cite{YanPZ17}.

Furthermore, people have recently realized that stochastic hybrid models have advantages to
capture discrete movements (such as random environment, market trends, interest rates, business
cycles, etc.) in the insurance market. For example, we can consider an insurance market with two
modes to represent the dynamic insurance cycle. Market mode 1 represents a ``soft"
market, where the investment return is high and the premium rate is low. While market mode
2 represents a ``hard" market, where the investment return is low and the premium rate is high. In different
markets, insurance companies adopt completely different strategies in investment and policy sales.
The insurance company is more likely to expand its business and write more policies when market is
in mode 1. The potential loss from the difference between premium income and claims
are compensated by the high investment returns.
The hybrid systems enable the consideration of
the coexistence of continuous dynamics and discrete
events in the systems. To reflect the hybrid feature, one of the recent trends is to use a finite
state Markov process to describe the transitions among different regimes. The Markov-modulated
switching systems are
known as regime-switching systems. The formulation of regime-switching
models is a more general and versatile framework to describe the complicated financial markets and their inherent uncertainty and randomness.
 Because the control strategies are affected by the asset prices on the stock market and economic trends change quickly,
Markovian regime-switching processes were introduced widely to capture movements of random environment.
In \cite{WYW10}, the optimal
dividend and proportional reinsurance strategy under utility
criteria are studied for the regime-switching compound Poisson
model. \cite{SotomayorC} studies the optimal dividend problem in the regime-switching model when the dividend rates are bounded, unbounded,
and when there are fixed costs and taxes corresponding to the dividend payments.
\cite{Zhu2014} studied the dividend optimization for a regime-switching diffusion model with restricted dividend rates. A comprehensive study of
switching diffusions with ``state-dependent'' switching is in
\cite{YinZ10}. \cite{BensoussanSYY} provided closed-form Nash equilibrium for a mixed regime-switching Cram\' er-Lundberg diffusion approximation process; see also related works
 \cite{BYY2012} and \cite{BHYY} for regime-switching models for real options and real  options with competition.
 \cite{JinYW13} designed numerical methods for a zero-sum stochastic differential reinsurance game with regime-switching.

In this work, we are concerned with an insurance market including two insurance companies. The two competing insurance companies adopt optimal investment and reinsurance strategies to manage the insurance portfolios. The surplus process of each insurance company is subject to the randomness of the market. Following the work of \cite{BensoussanSYY}, the randomness of the market is modelled by a continuous-time finite-state Markov chain and an independent market-index process. Nevertheless, we model the surplus process as a regime-switching jump-diffusion process, in lieu of a mixed regime-switching Cram\' er-Lundberg diffusion approximation process. This allows us to work with both proportional and excess-of-loss reinsurance policies. Equilibrium strategies are studied by solving
a system of HJI (Hamilton-Jacobi-Isaacs) equations for the value functions of various players, derived from the principle of dynamic programming. Owing to the inclusion of the random switching environment and jump processes, the system of HJI equations becomes more complicated and
 closed-form solutions are virtually impossible to obtain. Thus, we adopt the Markov chain approximation method (MCAM) developed in \cite{KushnerD}
 to deal with a system of HJI  (HJI) equations arising from the associated game problems.  The convergence of the approximation sequence to the jump process and the convergence of the value function will be
 established. In the actual computation, we will
  use our approximation schemes for constant absolute risk aversion (CARA) insurers.

As far as the significance of the contributions is concerned, this paper reveals clearly the advantage in the following aspects.
First, the problem of finding the optimal investment and reinsurance strategy between two insurance companies
under different insurance market models involves a stochastic hybrid system. Closed-form solutions can only be obtained in special cases in \cite{BensoussanSYY} in a diffusion approximation model.
Comparing with work in \cite{BensoussanSYY}, a more
versatile jump-diffusion
regime-switching non-zero-sum game model is formulated. \cite{JinYW13} only studies a zero-sum reinsurance game, where the state process is degenerated to one dimension. While in the
non-zero-sum game of the current work, a high-dimension system of HJI equations is obtained, which adds much difficulty to construct the approximating Markov chain. Second,
the
aforementioned literature
mainly considers
proportional reinsurance strategies in the reinsurance games. \cite{PunW16} derives
the Nash equilibrium of the robust non-zero-sum stochastic differential reinsurance game where the reinsurance strategies are proportional. Due to the unboundedness of the excess-of-loss strategies,
it is very hard to develop the Nash equilibrium strategies for an reinsurance game. In addition to the investment, finding the explicit optimal game strategies are virtually impossible. To the best of our
knowledge, this is the first paper to study the non-zero-sum investment and reinsurance games in a jump-diffusion regime-switching with nonproportional reinsurance strategies. Third, the reinsurance
premium principle, which is a generalized nonlinear functional in our model, adds further difficulties to design the game strategies. \cite{MengSY13} and \cite{MengLJ15} investigate the nonlinear
premium principles in a simplified reinsurance game formulation. Our models considers the investment and reinsurance game in more complex formulation and presents the impact of the nonlinearity
of the premium principle on investment and reinsurance strategies. Hence, the constructed numerical algorithm shows its advantage in finding the game strategies under our complex
stochastic hybrid system.

The rest of the paper is organized as follows. A generalized formulation for surplus processes and the
 associated control problem are presented in Section \ref{sec:form}. We design the numerical algorithm based on MCAM in Section \ref{sec:num}. A directly constructed approximating Markov chain is constructed and is proved to be locally consistent with the original processes. Section \ref{sec:con} deals with the convergence of the approximated process and the value functions. Numerical examples are
 reported in Section \ref{sec:exam} to illustrate the performance of the method. Section \ref{sec:rem} concludes the paper with a few more remarks.

\section{Formulation}\label{sec:form}
Let $(\Omega, \F, \F_t,P)$ be a complete filtered probability space,  where the filtration $\{ \F_t\}_{t \ge 0}$ satisfies the usual condition  and $\F_t = \F_t^\al \vee \F_t^W$.
We work with a finite horizon  $[0,T]$, where
$T< \infty$ is a positive
real number.
The processes $\{ \al(t) \}_{t \ge 0}$ and $\{ W(t) \}_{t \ge 0} =\{W_Z(t), W_S(t) \}_{t \ge 0}$ are defined below. We use $K$ as a generic constant throughout the paper, whose value may change for different appearances.

\subsection{Insurance Models}
We are considering an insurance market including two competing insurance companies. Each of the two insurance companies adopts optimal investment and reinsurance strategies to manage the insurance portfolios. The surplus process of each insurance company is subject to the random fluctuation of the market. Following the work of \cite{BensoussanSYY}, the randomness of the market is modelled by a continuous-time finite-state Markov chain and an independent market-index process.

To delineate the random economy environment and other random economic factors, we use a continuous-time Markov chain $\al(t)$ taking values in a finite space $\M=\{1,\ldots,m\}$. The states of economy are represented by the Markov chain $\al(t)$.
Let the continuous-time Markov chain $\al(t)$ be generated by $Q=(q_{ij})\in \rr^{m\times m}$. That is,
\beq{gen}
\pr \{ \al(t + \delta )=j | \al(t)=i, \al(s), s \leq t \}
= \left \{ \barray q_{ij}\delta + o(\delta),  \qquad  \hfill {\rm if}\ j \neq i, \\
1 + q_{ii}\delta + o(\delta), \qquad  \hfill {\rm if}\ j = i,
\earray \right.
\eeq
where $q_{ij} \geq 0$ for $i, j =1, 2, \ldots, m$ with $j \neq i$ and $\sum_{j \in \M}q_{ij} =0$
for each $i \in \M$.

Furthermore, we are considering the insurance portfolios in a financial market with a market index $Z(t)$, whose  prices satisfies
\beq{index}
dZ(t)=\mu_Z (t, Z(t))dt + \sigma_Z(t, Z(t))dW_Z(t),
\eeq
where $W_Z(t)$ is a standard Brownian motion.
 Denote by  $\{\F_t^{W_Z} \}$ the filtration generated by Brownian motion $\{ W_Z(t) \}_{t \ge 0}$.
 We note that $Z(t)$ captures the dynamics of the financial market. The cash flows of the insurance companies such as the premiums of insurance policies, claims, and expenses, are subject to the performance of financial market. Hence, the key parameters of the surplus process are defined as functionals of both the finite-state Markov chain $\al(t)$ and market index $Z(t)$.

Following the classical Cram\'er-Lundberg process, we assume that $\widehat X_k(t), k\in\{1,2\}$,  the surplus of insurance company $k$ without investment
and reinsurance satisfies
\beq{classic} d \widehat X_k(t) = c_k(\al(t), Z(t))dt - Y_k(t), \
\ t \geq 0, \eeq
where $\widehat X_k(0):=\hat x_k$ is the initial surplus,  $c_k(\al(t), Z(t))$ is the rate of premium, and $Y_k(t)=\disp\sum_{i=1}^{N_k(t)} A^k_i$ is a compound Poisson process with the claim size $A^k_i$ with $\{A^k_i: i>1\}$ being a sequence of positive, independent and identically distributed random variables.

In this paper, we consider a Poisson measure in lieu of the traditionally used Poisson process. Suppose $\Theta \subset\rr_+$ is a compact set and the function $q_k\cd$ is the magnitude of the claim sizes.
\beq{Ni-def}N_k(t, H)= \hbox{ number of claims on } [0,t]
 \hbox{ with claim size taking values in }
 H \in \Theta,  \eeq
  counts the number of claims up to time $t$, which is a Poisson counting process. For $k=1,2$, $Y_k(t)$ is jump processes representing claims for each company with arrival rate $\la_k$. Note that claim frequencies depends on the economy and financial market states. The function $q_k(\al(t), Z(t), \rho_k)$ is assumed to be the magnitude of the claim sizes, where $\rho_k$ has distribution $\Pi_k\cd$, and $q_k(i, \cdot,\rho_k)$ is continuous for each $\rho_k$ and each $i\in \M$. At different regimes and financial market states, taking into consideration of random environment, the values of  $q_k(i, \cdot,\rho_k)$ could be much different. Then the Poisson measure $N_k(\cdot, \cdot)$ has intensity $\la_k dt\times \Pi_k(d\rho_k)$ where $\Pi_k(d\rho_k)=f(\rho_k)d\rho_k$.

Let $\nu_{k,n}$ denote 
the time of the $n$-th claim and $\zeta_{k,n} = \nu_{k,n+1} -\nu_{k,n}$. Let $\{\zeta_{k,v}, \rho_k, v \geq n\}$ be independent of $\{X_k(s), \al(s), s \leq \nu_{k,n}, \zeta_{k,v}, \rho_k, v \leq n\}$. Then the $n$-th claim term is $q_k ( \al(\nu_{k,n}), Z(\nu_{k,n}), \rho_k )$, and the claim amount of $Y_k$ can be written as $\disp Y_k(t) = \sum_{\nu_{k,n} \leq t}  q_k(\al(\nu_{k,n}),(\nu_{k,n}), \rho_k)$.


\subsection{Reinsurance and Investment}
Let $a_k(t)$ be the $\F_t$-progressively measurable process valued in $[0,1]$, an exogenous retention level, which is a control chosen by the insurance company representing the reinsurance policy and $g(a_k)$ is the reinsurance premium rate. Denote by $\A_k := \{a_k(t): 0 \le a_k(t) \le 1, 0 \le t \le T\}$, the set of reinsurance strategies of insurer $k$.
Recall that $A^k_i$
is
the size of the $i$th claim. Let $A^k_i(a_k)$ be the fraction of each claim paid by the primary insurance company. Then the aggregation claim amount paid by the primary insurance company is denoted as $Y_k^{a_k}(t)$.

\begin{rem}
{\rm Note that both the claim frequencies and severities
are depending on the Markov regimes and market index. It is a more general formulation compared with the work in \cite{BensoussanSYY}, where only the claim frequencies  depends on the  Markov regimes and market index. Therefore, the diffusion approximation in
 \cite{BensoussanSYY} is simplified due to the state-independence of claim severity. Further, with the compound Poisson jumps, the surplus process forms a controlled jump-diffusion regime-switching process. We aim to find
  optimal reinsurance strategies under the
  jump-diffusion regime-switching process formulation numerically.
}
\end{rem}

The insurance companies invest in both risk-free assets $S_0(t)$ and risky assets $S(t)$ with prices satisfying
  \beq{fin}  \left \{ \barray \ad
  \frac{dS_0(t)}{S_0(t)}=r(\al(t)) dt, \\ \ad
  \frac{dS(t)}{S(t)}=\mu_S (\al(t), Z(t))dt + \sigma_S(\al(t), Z(t))dW_S(t),
  \earray \right. \eeq
  where $r(\al(t))$ and $\mu(\al(t), Z(t))$ are the return rates of the risk-free and risky assets, respectively; $\sigma_S(\al(t), Z(t))$ is the corresponding volatility; $W_S(t)$ is a standard Brownian motion independent of $W_Z(t)$. For $k=1,2,$ the investment   behavior of the insurer $k$ is modelled as a portfolio process $b_k(t)$, where $b_k(t)$ is invested in the risky asset $S(t)$. Let $\B _k = \{ b_k(t): 0 \leq t \leq T \}$ denote the set of investment strategies of insurer $k$.

Combining the reinsurance and investment strategies, the surplus process of the insurance company $k$, denoted by $ \wdt X_k(t)$, follows
\beq{sur}
\begin{cases}
  d\wdt X_k(t) = \left\{r(\al(t)) \wdt X_k(t)+b_k (t)[\mu_S(\al(t), Z(t))-r(\al(t))]+c_k(\al(t), Z(t))-g(a_k(t))\right\}dt \\
\hspace{45pt} + b_k (t) \sigma_S(\al(t), Z(t))dW_S(t)-dY_k^{a_k}(t), \\
 \wdt X_k(0)=\wdt x_k,
  \end{cases} \eeq
where
$ \disp Y_k^{a_k}(t)=\int_0^t \int_{\rr_+}\tilde q_k(q_k, a_k)N_k(dt,d\rho_k)$ and $\tilde{q}_k(q_k,a_k)$ is the magnitude of the claim sizes \wrt the surplus process.

In this work, we  model the competition of two insurance companies with investment and reinsurance schemes in finite time horizon using a game theoretic formulation. The performance of each company is measured by the relative performance of their surpluses against their competitor's. Thus, the competition between the two companies formulates a game with two players, each of which can adjust its reinsurance strategies based on the competitor's
scheme. Let the relative surplus performance
for insurance company $k$ be $ X_k(t):=\wdt X_k(t)-\kappa_k \wdt X_l(t)$. Hence, $ X_k(t)$
is governed by the following dynamic system
  \beq{re-sur} \barray \ad dX_k(t)\\
  \ad=\sum_{i\in {\cal M} }  I_{\{\al(t)=i\}} \bigg\{r(i) X_k(t)+(b_k (t)-\kappa_k b_l(t))[\mu_S(i, Z(t))-r(i)]+c_k(i, Z(t))-\kappa_k c_l(i, Z(t)) \\
   \aad \quad  -(g(a_k(t))-\kappa_k g(a_l(t))) dt
   + (b_k (t)-\kappa_k b_l(t)) \sigma_S(i, Z(t))dW_S(t)\bigg\}-dY_k^{a_k}(t)+\kappa_k dY_l^{a_l}(t),
  \earray \eeq
  where
$$ Y_k^{a_k}(t)=\int_0^t \int_{{\rr_+}}\tilde q_k(q_k, a_k)N_k(dt,d\rho_k).$$

 \subsection{Proportional Reinsurance}
 We allow the insurance companies to continuously reinsure a fraction of its claim with the retention level $a_k\in[0,1]$, $k=1,2$.
 Note that
 $a_k$ is the exogenous retention level, and the control chosen by the insurance company for the reinsurance policy. Then $\tilde q(q_k, a_k)=a_k(t) q_k(\al(t), Z(t), \rho_k)$.
  We have $$Y_k^{a_k}(t)=\disp\sum_{i=1}^{N_k(t)} A^k_i(a_k)=\sum_{i=1}^{N_k(t)} a_k A^k_i .$$

Considering the proportional reinsurance strategies, for $k=1, 2$. The relative surplus process of the insurance company $k$, under the reinsurance and investment, follows
 \beq{re-sur1}
 \begin{cases}
 &dX_k(t) =\bigg\{r(\al(t)) X_k(t)+(b_k (t)-\kappa_k b_l (t))[\mu_S(\al(t), Z(t))-r(\al(t))]+c_k(\al(t), Z(t))\\
& \quad   -\kappa_k c_l(\al(t), Z(t)) -[g(a_k(t))-\kappa_k g(a_l(t))]\bigg\}dt
   + (b_k (t)-\kappa_k b_l(t)) \sigma_S(\al(t), Z(t))dW_S(t)\\
 & \quad -a_k(t)\int_{\rr_+}q_k(\al(t),Z(t),\rho_k)N_t(dt,d\rho_k) +\kappa_k a_l(t)\int_{\rr_+}q_l(\al(t),Z(t),\rho_l)N_l(dt,d\rho_l) \\
 X_k(0)&= \wdt x_k-\kappa_k  \wdt x_l.
\end{cases}  \eeq

 \subsection{Excess-of-loss Reinsurance}
We allow the insurance companies to continuously reinsure its claim and pay all of the claims up to a pre-given level of amount (termed retention level). We still let $a_k$, $k=1, 2$ be the retention level chosen by the insurance company to determine the reinsurance policy. We have that $$Y_k^{a_k}(t)=\disp\sum_{i=1}^{N_k(t)} A^k_i(a_k)=\sum_{i=1}^{N_k(t)} (A^k_i \wedge a_k).$$
Then $\tilde q_k(q_k, a_k)=q_k(\al(t),Z(t),\rho_k)\wedge a_k(t)$.

Considering the excess-of-loss reinsurance strategies, for $k=1, 2$, the relative surplus process of the insurance company $k$, under the reinsurance control and investment, follows
  \beq{re-sur2} \begin{cases}
  dX_k(t)&=\bigg\{r(\al(t))X_k(t)+(b_k (t)-\kappa_k b_l(t))[\mu_S(\al(t), Z(t))-r(\al(t))]+c_k(\al(t), Z(t)) \\
&\quad -\kappa_k c_l(\al(t), Z(t)) -[g(a_k(t))-\kappa_k g(a_l(t))]\bigg\}dt
   + (b_k (t)-\kappa_k b_l(t)) \sigma_S(\al(t), Z(t))dW_S(t)\\
  &\quad  -\int_{\rr_+}(q_k(\al(t),Z(t),\rho_k)\wedge a_k)N_k(dt,d\rho_k) +\kappa_k\int_{\rr_+}(q_l(\al(t),Z(t),\rho_l)\wedge a_l)N_l(dt,d\rho_l) \\
 X_k(0)& =\wdt x_k-\kappa_k \wdt  x_l.
  \end{cases} \eeq

  \subsection{Control Problem}
  For $k=1,2$,
insurer $k$ has a utility function $U_k: \rr \ar \rr$, where $U_k$ is assumed to be increasing, strictly concave, and satisfies Inada conditions, i.e.,
$$\partial_x U_k(-\infty) = +\infty, \qquad \partial_x U_k(+ \infty) =0.$$
Following the work  \cite{EspinosaT}, the insurer $k$ aims to maximize the expected utility of his relative performance at the terminal time $T $ by adopting a pair of investment and reinsurance strategy $u_k=(a_k, b_k) \in \A_k \times \B_k$, denote $\U_k \triangleq \A_k \times \B_k$. For an arbitrary pair of admissible control $u =(u_1, u_2) \in \U \triangleq\U_1 \times \U_2$, the objective function is
 \beq{objective} \barray
 J^k (t,x_k,z,i,u) \ad =
 \ex \left[ U_k \left( (1-\kappa_k)\wdt X_k(T)+\kappa_k(\wdt X_k(T)-\wdt X_l(T))\right)\right]
 \\ \ad = \ex \left[ U_k\left( \wdt X_k(T)-\kappa_k \wdt X_l(T)\right)\right], \earray \eeq
for $k\neq l \in \{1, 2\}$. For $k=1,2$, $\kappa_k$ measures the sensitivity of insurer $k$ to the performance of his competitor.

The control $u_k=(a_k,b_k)$ with $k\in\{1, 2\}$ is said to be {\em admissible} if $a_k$ and $b_k$ satisfy
 \begin{itemize}
 \setlength{\baselineskip}{0.14in}
\parskip=2pt
 \item[(i)] $a_k(t), b_k(t)$ are nonnegative for any $t\ge 0$,
 \item[(ii)] Both $a_k$, $b_k$ are adapted to ${\cal F}_t$.
 \item[(iii)] $J^k (t,x_k,z,i,u) < \infty$ for any admissible pair $u_k=(a_k, b_k)$.
 \end{itemize}

 For $k =1,2,$ let $\bb(\U_k \times [0,\infty))$ be the $\sg-$ algebra of Borel subsets of $\U_k \times [0, \infty)$.
 We use a relaxed control formulation; see \cite{KushnerD} for a definition and more discussions. Recall that
 an \textit{admissible relaxed control} $m_k \cd$ is a measure on $\bb(\U_k \times [0, \infty))$ such that $m_k(\U_k \times [0,t)) =t$ for each $t \ge 0$. With the given probability space, we say that $m_k \cd$ is an admissible relaxed (stochastic) control for $\F$, if $m_k(\cdot, \omega)$  is a deterministic relaxed control with probability one and if $m_k(A \times [0,t])$ is $\F_t$-adapted for all $A \in \bb(\U_k)$.

 Given a relaxed control
 $m_k \cd$ of $u_k \cd$, we define the derivative $m_{t,k} \cd$ such that
 $$m_k(K) = \int_{\U_k \times [0,\infty)} I_{\{ (u_k, t) \in K \}} m_{t,k}(d\phi_k)dt$$
 for all $K \in \bb(\U_k \times [0, \infty))$, and that for each $t$, $m_{t,k} \cd$ is a measure on $\bb(\U_k)$ satisfying $m_{t,k} (\U_k) =1$. For example, we can define $m_{t,k} \cd$ in any convenient way for $t =0$ and as the left-hand derivative for $t >0$,
 $$m_{t,k}(A) = \lim_{\varrho \to 0} \dfrac{m_k(A \times [t -\varrho,t])}{\varrho}, \qquad \forall A \in \bb(\U_k).$$
 Note that $m_k (d \phi_k dt) = m_{t,k}(d\phi_k) dt$. It is natural to define the relaxed control representation $m_k \cd$ of $u_k \cd$ by
 $m_{t,k}(A) = I_{\{u_k(t) \in A \}}$  $\forall A \in \bb(\U_k).$
 Define the relaxed control $m \cd = (m_1 \cd \times m_2 \cd)$ with derivative $m_t \cd = m_{t,1}\cd \times m_{t,2} \cd$. Thus $m \cd$ is a measure on the Borel sets of $\U \times [0, \infty)$.

  \subsection{Nash Equilibrium}
 A Nash equilibrium $u^*=(u_1^*, u_2^*)\in \U $ is achieved such that
\beq{nash} \barray \ad
\ex[U_1(\wdt X_1^{u_1}(T)-\kappa_1 \wdt X_2^{u_2^*}(T))] \leq \ex[U_1(\wdt X_1^{u_1^*}(T)-\kappa_1 \wdt X_2^{u_2^*}(T))], \\ \ad
\ex[U_2(\wdt X_2^{u_2}(T)-\kappa_2 \wdt X_1^{u_1^*}(T))] \leq \ex[U_2(\wdt X_2^{u_2^*}(T)-\kappa_1 \wdt X_1^{u_1^*}(T))].
\earray \eeq
For $\al(t)=i \in \M$, $Z(t)= z$, and $X_k(t)=x_k$, where $0\leq t \leq T$ and $k\neq l \in \{1, 2\}$, the value function of insurance company $k$ follows
\beq{value}
V^k(t, x_k, z, i)=\sup_{u_k \in \H_k} \ex[U_k(\wdt X_k^{u_k}(T)-\kappa_k \wdt X_l^{u_l^*}(T))],
\eeq
where $V^k(\cdot, \cdot, \cdot, \cdot)$ is the value function in $\rr_+ \times \rr \times \rr \times \M $.

To obtain the system of Hamilton-Jacobi-Bellman (HJB) equations, we assume the existence of optimal control. For an arbitrary $u_k\in \U_k$, $\al(t)=i \in\M$,  $k\neq l \in \{1, 2\},$ and
$V^k(\cdot,\cdot,\cdot, i)\in C^{2}(\rr_+ \times \rr \times \rr\times \M)$,
   define an integro-differential operator ${\cal L}^{u_k, u_l}$ by
   \beq{lth}\barray
  \Ll^{u_k, u_l}V^k(t,x_k,z,i)\ad= V^k_{x_k }(t,x_k,z,i)\bigg\{r(i)x_k+(b_k -\kappa_k b_l)[\mu(i, z)-r(i)]+c_k(i, z)-\kappa_k c_l(i, z)
   \\ \aad \quad    -(g(a_k)-\kappa_k g(a_l))\bigg\}+ V^k_{z}(t,x_k,z,i)m(t,z)
   \\ \aad \quad  +\frac{1}{2} (b_k -\kappa_k b_l)^2 \sigma^2_S(i, z)V_{x_k x_k}(t,x_k,z,i)+\frac{1}{2}\sigma^2_Z(t, z)V_{zz}(t,x_k,z,i)
   \\ \aad\quad -\la_k\int_{\rr_+} [V^k(t,x_k-\tilde q(i,z,\rho_k),z,i)- V^k(t,x_k,z,i)]f(\rho_k)d\rho_k
   \\ \aad\quad +\la_l\int_{\rr_+} [V^k(t,x_k+\kappa_k \tilde q(i,z,\rho_l),z,i)- V^k(t,x_k,z,i)]f(\rho_l)d\rho_l
   \\ \aad\quad +Q V(t,x_k,z,\cdot)(i) ,
   \earray\eeq
   where
   $$ Q V(t,x_k,z,\cdot)(i)=\sum_{j\neq i}q_{ij}(V(t,x_k,z,j)-V(t,x_k,z,i)).$$
  Formally, for $k\neq l \in \{1, 2\}$, we conclude that $V^k$ satisfies the following system of integro-differential HJI
  (Hamilton-Jacobi-Isaacs)
  equations: for each $ i\in\M$,
  \beq{hji}
\begin{cases}
  &V^k_t(t,x_k,z,i)+\disp \sup_{u_k\in \H_k}\Ll^{u_k, u_l^*}V^k(t,x_k,z,i)=0,
   \\ &  V^k(T,x_k,z,i)=U_k(x_k).
\end{cases}
 \eeq

 \section{Numerical Algorithm}\label{sec:num}
 We begin by construction a discrete-time, finite-state, controlled Markov chain to approximate the controlled diffusion process with regime-switching in the absence of jumps with the dynamic system
  \beq{nojump} \barray dX_k(t)
  \ad=\sum_{i\in \M }  I_{\{\al(t)=i\}} \bigg\{r(i) X_k(t)+(b_k (t)-\kappa_k b_l(t))[\mu_S(i, Z(t))-r(i)]+c_k(i, Z(t))-\kappa_k c_l(i, Z(t)) \\
   \aad \quad  -(g(a_k)+\kappa_k g(a_l)) dt
   + (b_k (t)-\kappa_k b_l(t)) \sg_S(i, Z(t))dW_S(t)\bigg\},\\
   dZ(t)\ad =\mu_Z (t, Z(t))dt + \sg_Z(t, Z(t))dW_Z(t)\\
  X_k(0) \ad = \wdt x_k - \kappa_k \wdt x_l \\
  Z(0)\ad = z_0.
  \earray \eeq

Because the value function depends on both the state $x$ and the time variable $t$, two stepsizes are needed.
That is, we need to discretize both the state and time.
We use $h >0$ as the stepsize of the state and $\dl >0$ as the stepsize for the time. In fact, for any given $T>0$, we use $N = N(\dl) = \lfloor T/\dl \rfloor$.

Let $e_i$ denote the unit vector in the $i$-th coordinate direction and $\rr^3_h$ denote the uniform $h$-grid on $\rr^3$; i.e. $\rr^3_h = \{ (x_1,x_2,z): (x_1,x_2,z) = h(k_1e_1+ k_2 e_2+k_3 e_3); k_1, k_2,k_3 =0, \pm 1, \pm 2,\dots \}$. We use $S_h = \rr^3_h$, denote $x = (x_1, x_2, z)$ and $y = (y_1, y_2, z^\ast)$.

We can rewrite the system in the short form as the following:
\beq{m-form} \barray
dX(t) = \mu(X(t), \al(t), u(t)) dt + \sg(X(t), \al(t), u(t)) dW(t),
\earray
\eeq
where
\bea
\ad\!\!\!\! \mu(X(t), \al(t), u(t)) \\
\ad =\begin{pmatrix}
 r(\al(t)) X_1(t)+(b_1 (t)-\kappa_1 b_2(t))[\mu_S(\al(t), Z(t))-r(\al(t))]+c_1(\al(t), Z(t))-\kappa_1 c_2(\al(t), Z(t))\\
  -(g(a_1(t))+\kappa_1 g(a_2(t))) \\
 r(\al(t)) X_2(t)+(b_2 (t)-\kappa_2 b_1(t))[\mu_S(\al(t), Z(t))-r(\al(t))]+c_2(\al(t), Z(t))-\kappa_2 c_1(\al(t), Z(t))\\
  -(g(a_2(t))+\kappa_2 g(a_1(t)))\\
\mu_Z (t, Z(t))
\end{pmatrix} \\
\ad\!\!\!\! \sg(X(t), \al(t), u(t))
\\
\ad = \diag \bigg((b_1 (t)-\kappa_1 b_2(t)) \sg_S(\al(t), Z(t)), (b_2 (t)-\kappa_2 b_1(t)) \sg_S(\al(t), Z(t)), \sg_Z(t, Z(t)) \bigg)\\
\ad\!\!\! W(t) = \bigg(W_S(t), W_S(t), W_Z(t) \bigg)',
\eea
where $A'$ is the transpose of $A$.
Let $ \{( \xi^{h,\dl}_n, \al^{h,\dl}_n , n < \infty \}$ be a controlled discrete-time Markov chain on $\rr^3_h \times \M$ and denote by $p^{h,\dl}_D \big((x,i), (y,j),n\dl |\phi \big)$  the transition probability from a state $(x,i)$ to another state $(y,j)$, for $\phi \in \U$. We use $u^{h, \dl}_n$ to denote the random variable that is the control action for the chain at discrete time $n$ and $p^{h,\dl}_D$ is so defined that the constructed Markov chain's evolution well approximates the local behavior of the controlled regime-switching diffusion \eqref{re-sur}.

For each $k =1,2$, we construct the transition probability $p^{h,\dl}_{k,D}\big((x,i),(y,j),n \dl | \phi \big)$ which is associated with $J^k(t,x,i,\phi) = \ex \left[ U_k\left(X_k(T)\right)\right]$ satisfying the followings:
  \beq{hji1}
\begin{cases}
  &J^k_t(t,x,i,\phi)+\disp \Ll^{u_k, u_l}J^k(t,x,i,\phi)=0,
   \\ &  J^k(T,x,i,\phi)=U_k(x_k).
\end{cases}
 \eeq
To figure out the form of $p^{h,\dl}_{k,D}\big( (x,i),(y,j)| \phi \big)$, we define a finite difference approximation to \eqref{hji1} as
\bea
& J^k_t (t,x,i,\phi) \ar \dfrac{J^k(t+\dl, x,i,\phi)-J^k(t,x,i,\phi)}{\dl},\\
& J^k_{x_k} (t,x,i,\phi) \ar \dfrac{J^k(t,x+he_k,i,\phi)-J^k(t,x,i,\phi)}{h} \\
&\quad \hbox{if}\  r(i)x_k+(b_k -\kappa_k b_l)[\mu_S(i, z)-r(i)]+c_k(i, z)-\kappa_k c_l(i, z) -(g(a_k)-\kappa_k g(a_l)) >0,\\
& J^k_{x_k} (t,x,i,\phi) \ar \dfrac{J^k(t,x,i,\phi)-J^k(t,x-he_k,i,\phi)}{h} \\
&\quad \hbox{if}\  r(i)x_k+(b_k -\kappa_k b_l)[\mu_S(i, z)-r(i)]+c_k(i, z)-\kappa_k c_l(i, z) -(g(a_k)-\kappa_k g(a_l)) <0, \\
& J^k_{x_k x_k} (t,x,i,\phi) \ar \dfrac{J^k(t,x+he_k,i,\phi) - 2J^k(t,x,i,\phi) + J^k(t,x-he_k,i,\phi) }{h^2} ,\\
& J^k_z (t,x,i,\phi) \ar \dfrac{J^k(t,x+he_3,i,\phi) - J^k(t,x,i,\phi)}{h} \quad \hbox{if}\quad  \mu_Z(t,z)>0, \\
& J^k_z (t,x,i,\phi) \ar \dfrac{J^k(t,x,i,\phi) - J^k(t,x-h e_3,i,\phi)}{h} \quad \hbox{if}\quad \mu_Z(t,z)<0, \\
& J^k_{zz}(t,x,i,\phi) \ar \dfrac{J^k(t,x+h e_3, i,\phi) -2J^k(t,x,i,\phi) + J^k(t,x-h e_3,i,\phi)}{h^2}.
\eea
To proceed, define
\beq{prob1}
\barray
& p^{h, \dl}_{k,D} \big( (x,i),(x \pm h e_k,i), n \dl | \phi \big) = \frac{\dl}{h} \bigg\{ r(i)x_k+(b_k -\kappa_k b_l)[\mu_S(i, z)-r(i)]+c_k(i, z)-\kappa_k c_l(i, z)\\
&\quad -(g(a_k)-\kappa_k g(a_l)) \bigg\}^\pm + \frac{\dl}{2h^2} (b_k -\kappa_k b_l)^2 \sigma^2_S(i,z), \\
& p^{h,\dl}_{k,D} \big( (x,i), (x \pm h e_3,i), n \dl | \phi \big) = \frac{\dl}{h} \mu_Z(t,z)^\pm + \frac{\dl}{2h^2} \sigma_Z^2(t,z), \\
& p^{h, \dl}_{k,D} \big( (x,i), (x ,j), n \dl | \phi \big) = q_{ij} \dl,\\
& p^{h, \dl}_{k,D} \big( (x,i), (x ,i), n \dl | \phi \big) = 1 + q_{ii} \dl - \frac{\dl}{h} \left| r(i)x_k+(b_k -\kappa_k b_l)[\mu_S(i, z)-r(i)]+c_k(i, z)-\kappa_k c_l(i, z) \right.\\
&\left.\quad -(g(a_k)-\kappa_k g(a_l)) \right| -\frac{\dl}{h} \mu_Z(t,z)-\frac{\dl}{h^2} (b_k -\kappa_k b_l)^2 \sigma^2_S(i,z)-\frac{\dl}{h^2} \sigma_Z^2(t,z), \\
& p^{h, \dl}_{k,D} \big( \cdot \big) = 0, \quad \hbox{otherwise},
\earray
\eeq
where $K^+ = \max\{K,0 \}$ and $K^- = \min \{ -K,0\}$.
By choosing $\dl$ and $h$ appropriately, we can have $p^{h,\dl}_{k,D} \big( (x,i), (x,i), n \dl | \phi \big)$ given in \eqref{prob1} nonnegative. Thus, $p^{h, \dl} (\cdot| \phi)$ are well-defined transition probability.

Next, we need to approximate the Poisson jumps for ensuring the local properties of claims for \eqref{re-sur}. We can rewrite the system in the matrix form as follow
\beq{mform}
d X(t) = \mu(X(t), \al(t), u(t)) dt + \sg(X(t), \al(t), u(t)) dB(t) + \wdt Y_1(t) e_1 +\wdt Y_2(t) e_2,
\eeq
where $\wdt Y_k(t)$ is the jump process w.r.t the surplus process $X_k(t)$, for $k=1,2$.

 The relative surplus process $X_k(t)$ is determined by two jump terms with the arriving rate $\la_k$ and $\la_l$, respectively. Denote by $R_k(t)$ the difference of the two jumps. That is,
$$R_k(t) = \int_{\rr_+} \wdt q_k (q_k, a_k) N_k(dt,d \rho_k) - \kappa_k\int_{\rr_+} \wdt q_l (q_l, a_l) N_l(dt, d \rho_l).$$
Since the difference of two Poisson processes is again a Poisson process, events in the new process $R_k(t)$ will occur according to a Poisson process with the rate $\la$. with $ \la= \la_k + \la_l$; and each event, independently, will be from the first jump process with probability $\la_k/(\la_k + \la_l)$, yielding the generic claim size
\bea
\wdt A^k = \begin{cases}
A^k(a_k), & \quad \hbox{with probability } \frac{\la_k}{\la_k+\la_l},\\\\
- \kappa_k A^l(a_l),& \quad \hbox{with probability } \frac{\la_l}{\la_k+\la_l}.
\end{cases}
\eea

Suppose that the current state is $\xi^{h,\dl}_n =x$, $\al^{h,\dl}_n =i$, and control is $u^{h,\dl}_n = \phi$. The next interpolation interval is determined by \eqref{prob1} and $\wdt q_{k}^h(q_k(i,z,\rho_k), a_k)$ is the nearest value of $\wdt q_k (q_k(i,z, \rho_k), a_k)$ so that $\xi^{h,\dl}_{n+1} \in S_h$. Then $|\wdt q_{k}^h(q_k(i,z,\rho_k), a_k) -\wdt q_k (q_k(i,z, \rho_k), a_k)| \ar 0$ as $h \ar 0$, uniformly in x. To present the claim terms, we determine the next case $(\xi^{h,\dl}_{n+1}, \al^{h,\dl}_{n+1})$ by noting:
\begin{enumerate}
\item No claims occur in $[n \dl, n \dl+\dl)$ with probability $1 - \la \dl + o(\dl)$, we determine $(\xi^{h,\dl}_{n+1}, \al^{h,\dl}_{n+1})$ by transition probability $p^{h,\dl}_{k,D} \cd$ as in \eqref{prob1}.
\item There is a claim of the relative surplus process $X_k(t)$ in $[n \dl, n \dl+\dl)$ with probability $\la \dl + o(\dl)$, we determine $(\xi^{h,\dl}_{n+1}, \al^{h,\dl}_{n+1})$ by
$$\xi^{h,\dl}_{n+1} = \xi^{h,\dl}_n -\wdt q_{k}^h(q_k(i,z,\rho_k), a_k)e_k, \qquad \al^{h,\dl}_{n+1} =\al^{h,\dl}_n.$$
\end{enumerate}
So, we define
\beq{pk}
p^{h,\dl}_k \big((x,i), (y,j)| \phi \big) = (1 -\la \dl + o(\dl)) p^{h,\dl}_{k,D}\big((x,i), (y,j)| \phi \big)+ (\la \dl + o(\dl)) \Pi_k\{\rho: \wdt q_h(i,z,\rho) e_k =x-y  \}.
\eeq

\begin{defn}
 {\rm A controlled Markov chain $\{(\xi^{h,\dl}_n, \al^{h,\dl}_n), n < \infty \}$ with the one-step transition probability $p^{h,\dl} \big((x,i), (y,j)| \phi \big)$ is given by
\beq{prob2}
p^{h,\dl} \big((x,i), (y,j)| \phi \big) = p^{h,\dl}_1 \big((x,i), (y,j)| \phi \big) I_{\{k=1 \}}+p^{h,\dl}_2 \big((x,i), (y,j)| \phi \big) I_{\{k=2 \}},
\eeq
where  $k$ is the index of the cost function.
}
\end{defn}

The piecewise constant interpolations $\xi^{h, \dl} \cd, \al^{h,\dl} \cd$, and $u^{h,\dl} \cd$ are defined as
\beq{ipl}
\xi^{h,\dl} (t) = \xi^{h, \dl}_n,\quad \al^{h,\dl}(t) = \al^{h,\dl}_n, \quad u^{h,\dl} (t) = u^{h,\dl}_n \quad \hbox{for } t \in [n \dl, n \dl+\dl)
\eeq
Use $\exh, \vah$, and $\prh$ to denote the conditional expectation, variance, and marginal probability given $\{\xi^{h,\dl}_\iota, \al^{\dl}_\iota,u^{h,\dl }_\iota, \iota \le n, \xi^{h,\dl}_n = x, \al^{h,\dl}_n =i, u^{h,\dl}_n = \phi \}$, respectively. Define the difference $\Dl \xi^{h, \dl}_n = \xi^{h, \dl}_{n+1} -\xi^{h, \dl}_n$.

With the approximation of the Markov chain constructed above, we can obtain an approximation of the utility function as follows:
\beq{dp}
\barray
J^{k, \dl}(n\dl, x, i, u^{h,\dl}) &= \disp(1- \la \dl + o(\dl)) \sum_{(y,j)} p^{h,\dl}_k\bigg((x,i),(y,j)|u^{h,\dl} \bigg) J^{k,\dl}(n\dl+\dl, y,j,u^{h,\dl})\\
&+ \disp(\la_k \dl+ o(\dl)) \int_{\rr_+}J^{k,\dl}(n \dl+\dl, x_k -\wdt q_{k}^h(q_k(i,z,\rho_k), a^{h,\dl}_k), z, i, u^{h, \dl}) \Pi_k(d \rho_k)\\
&+ \disp (\la_l \dl + o(\dl)) \int_{\rr_+}J^{k,\dl} (n \dl + \dl, x_k + \kappa_k \wdt q_l^h(q_l(i,z,\rho_l), a^{h,\dl}_l),z,i,u^{h,\dl}) \Pi_l(d \rho_l).
\earray
\eeq
Moreover,
\beq{val}
\barray
V^{k,h,\dl}(n \dl, x,i)=\disp\sup_{u_{k,n}^{h,\dl}}J^{k,h,\dl}(n\dl, x, i, u^{h,\dl}).
\earray
\eeq

\begin{defn} {\rm
The sequence $\{(\xi^{h,\dl}_n, \al^{h,\dl}_n) \}$ is said to be {\bf locally consistent}, if it
\begin{enumerate}
\item There is a transition probability $p^{h,\dl}_D$ is {\bf locally consistent} in the sense
\beq{cst}
\barray
&\exh [\Dl \xi^{h, \dl}_n ]= \mu^{h,\dl}(x,i,\phi)  \dl + o(\dl), \\
&\vah [\Dl \xi^{h, \dl}_n]=  \sg^{h,\dl} (x,i,\phi)\dl +o(\dl).
\earray
\eeq
where
\bea
&\mu^{h,\dl}(x,i,\phi) =\begin{pmatrix}
r(i) x_1+(b_1 -\kappa_1 b_2)[\mu_S(i, z)-r(i)]+c_1(i, z)-\kappa_1 c_2(i, z)  -g(a_1)+\kappa_1 g(a_2)\\
r(i) x_2 +(b_2 -\kappa_2 b_1)[\mu_S(i,z) -r(i)]+c_2(i,z)- \kappa_2 c_1(i,z) -g(a_2) + \kappa_2 g(a_1)\\
\mu_Z(t,z)
\end{pmatrix}\\
& \sg^{h,\dl} (x,i,\phi) = \diag \bigg((b_1 -\kappa_1 b_2)^2 \sg_S^2(i,z), (b_2 -\kappa_2 b_1)^2 \sg_S^2(i,z), \sg_Z^2(t,z) \bigg).
\eea
\item The one-step transition probability $p^{h,\dl} ((x,i),(y,j) |\phi)$ for the chain can be represented in the factored form:
\bea
p^{h,\dl}((x,i),(y,j)|\phi) = (1- \la\dl +o(\dl)) p^{h,\dl}_D((x,i),(y,j)|\phi) + (\la \dl +o(\dl)) \Pi \{ \rho: q_h (i,z,\rho) = x-y\}
\eea
\end{enumerate}
}
\end{defn}

\section{Convergence of Numerical Approximation}\label{sec:con}
\subsection{Representations of Approximation Sequences}
To proceed, we first show that the constructed Markov chain is locally consistent. This ensures that our approximation is reasonable in certain sense.

\begin{lem} The Markov chain $\xi^{ h,\dl}_n$ with transition probabilities $p^{h, \dl} \cd$ defined in \eqref{prob2} is locally consistent with the stochastic differential equation in \eqref{re-sur}.
\end{lem}

\para{Proof.}
Define $p^{h,\dl}_D \cd = p^{h,\dl}_{1,D} \cd I_{\{ k=1\}} + p^{h,\dl}_{2,D} \cd I_{\{k=2\}}$Using \eqref{prob1} and \eqref{prob2}. It is easy to see that
\bea
\exh [\Dl \xi^{h, \dl}_n] \ad =  \exh [\Dl \xi^{ h,\dl}_n I_{\{ k =1\}}] + \exh [\Dl \xi^{ h,\dl} I_{\{ k =2\}}]\\\\
\ad= \begin{pmatrix}
r(i) x_1+(b_1 -\kappa_1 b_2)[\mu_S(i, z)-r(i)]+c_1(i, z)-\kappa_1 c_2(i, z)  -g(a_1)+\kappa_1 g(a_2)\\
r(i) x_2 +(b_2 -\kappa_2 b_1)[\mu_S(i,z) -r(i)]+c_2(i,z)- \kappa_2 c_1(i,z) -g(a_2) + \kappa_2 g(a_1) \\
\mu_Z(t,z)
\end{pmatrix} \dl.
\eea
Likewise, we obtain $\exh \big[\Dl \xi^{h,\dl}_n (\Dl \xi^{h,\dl}_n)' \big]$ and $\vah [\Dl \xi^{h, \dl}_n]$. \qed

Let $\xi^{h,\dl}(0)=x, \al^{h,\dl}(0) = i$. Define the relaxed control representation $m^{h,\dl}_k \cd$ of $u^{h,\dl}_k \cd$ by using its derivative $m^{h,\dl}_{t,k} (A) = I_{\{ u^{h,\dl} (t) \in A \}}$.
 Let $H^{h,\dl}_n$ denote the event that $\xi^{h,\dl}_n, \al^{h,\dl}_n$ is determined by the case of ``no claim occurs" and use $T^{\dl}_n$ to denote the event of ``one claim occurs". Let $I_{H^{h,\dl}_n}$ and $I_{T^{h,\dl}_n}$ be corresponding indicator functions, respectively. Then $I_{H^{h,\dl}_n}+I_{T^{h,\dl}_n} =1 $ and we can write
\beq{mc1}\barray
\xi^{h,\dl} (t) &\disp =\xi^{h,\dl} (0)  + \sum_{r=0}^{\lf t/\dl \rf -1} [\Dl \xi^{h,\dl}_r I_{H^{h,\dl}_r} + \Dl \xi^{h,\dl}_r I_{T^{h,\dl}_r}]\\
&\disp = x + \sum_{r=0}^{\lf t/\dl \rf -1} \exr[\Dl \xi^{h,\dl}_r I_{H^{h,\dl}_r}] + \sum_{r=0}^{\lf t/\dl \rf -1} (\Dl\xi^{h,\dl}_r -\exr \Dl \xi^{h,\dl}_r) I_{H^{h,\dl}_r} + \sum_{r=0}^{\lf t/\dl \rf -1} \Dl \xi^{h,\dl}_r I_{T^{h,\dl}_r}.
\earray\eeq
Define $\F^{h, \dl}_n$ as the smallest $\sg$-algebra generated by $\{ \xi^{ h,\dl}_r, \al^{h, \dl}_r, m^{h, \dl}_r, H^{h,\dl}_r, r \leq n \}$ and $\F^{h,\dl}_t$ as the smallest $\sg$- algebra generated by $\{ \xi^{h, \dl} (s), \al^{h,\dl} (s), m^{h, \dl}(s), H^{h, \dl} (s), s \leq t  \}$.

 For  $k,l =1,2$ and  $k \neq l$,
denote
\beq{eq:not}
\barray
  M^{h,\dl} (t)\ad =\disp\sum_{r=0}^{\lf t/\dl \rf -1} (\Dl\xi^{h,\dl}_r -\exr \Dl \xi^{h,\dl}_r) I_{H^{h,\dl}_r} \\
  Y^{h,\dl}_k (t)\ad = -\disp\sum_{r=0}^{\lf t/\dl \rf-1} [\Dl \xi^{h,\dl}_r]' e_k I_{T^{h,\dl}_r}\\
 \ad = \int_0^t \int_{\rr_+} \wdt q^h_k (q_k( \al^{h,\dl}(s),Z^{h,\dl}(s), \rho_k), a^{h,\dl}_k(s))  N^{h,\dl}_k(ds, d\rho_k)\\
\aad\quad -  \int_0^t \int_{\rr_+} \ka_k \wdt q^h_l(q_l(\al^{h,\dl}(s), Z^{h,\dl}(s), \rho_l), a^{h,\dl}_l(s)) N^{h,\dl}_l(ds,d\rho_l).
\earray
\eeq
Then $M^{h,\dl}(t)$ is a martingale with respect to $\F^{h,\dl}_{\lf t/\dl \rf}$. Now, we represent $M^{h,\dl}(t)$ similar to the diffusion term in \eqref{re-sur}. Define $W^{h,\dl}_k \cd$ as
\beq{BM} \barray
W^{h,\dl}(t)  &=\disp \sum_{r=0}^{\lf t/\dl \rf -1} [ \sg^{h,\dl} (x,\phi, i) ]^{-1}(\Dl\xi^{h,\dl}_r -\exr \Dl \xi^{h,\dl}_r) I_{H^{h,\dl}_r}\\
&= \disp\int_0^t  [ \sg^{h,\dl}  (\xi^{h,\dl}(s), \al^{h,\dl}(s), u^{h,\dl}(s) )]^{-1}  dM^{h,\dl}(s) .
\earray \eeq

The local consistency leads to
\bea
&\disp \sum_{r=0}^{\lf t/\dl \rf -1} \exr [\Dl \xi^{h,\dl}_r I_{H^{h,\dl}_r}]\\
& = \disp \sum_{r=0}^{\lf t/\dl \rf -1}\begin{pmatrix}
r(\al^{h,\dl}_r) \xi^\dl_r e_1+(b^{h,\dl}_{1,r} -\ka_1 b^{h,\dl}_{2,r})[\mu(\al^{h,\dl}_r, \xi^{h,\dl}_r e_3)-r(\al^{h,\dl}_r)]\\
+c_1(\al^{h,\dl}_r, \xi^{h,\dl}_r e_3)-\ka_1 c_2(\al^{h,\dl}_r, \xi^{h,\dl}_r e_3)  -g(a^{h,\dl}_{1,r})+\ka_1 g(a^{h,\dl}_{2,r})\\
r(\al^{h,\dl}_r)\xi^{h,\dl}_r e_2 +(b^{h,\dl}_{2,r} -\ka_2 b^{h,\dl}_{1,r})[\mu(\al^{h,\dl}_r, \xi^{h,\dl}_r e_3) -r(\al^{h,\dl}_r)]\\
+c_2(\al^{h,\dl}_r,\xi^{h,\dl}_r e_3)- \ka_2 c_1(\al^{h,\dl}_r,\xi^{h,\dl}_r e_3) -g(a^{h,\dl}_{2,r}) + \ka_2 g(a^{h,\dl}_{1,r})\\
\mu_Z(r \dl,\xi^{h,\dl}_r e_3)
\end{pmatrix}  \dl I_{H^{h,\dl}_r} + o(\dl) I_{H^{h,\dl}_r}\\
& = \disp\int_0^t \mu^{h,\dl} (\xi^{h,\dl}(s), \al^{h,\dl}(s), u^{\dl}(s))  ds + \e^{h,\dl}(t).
\eea
For each $t$,
$
\ex [\hbox{number of } r : \nu^{h,\dl}_{k,r} \leq t ] = \la_k t
$
as $h,\dl \ar 0$. This implies that we can drop $I_{H^{h,\dl}_r}$ with no effect on the above limit.

As a consequence, we can rewrite \eqref{mc1} as following:
\beq{mc2} \barray
\xi^{h,\dl}(t)
& = \disp x + \int_{0}^{t} \int_\U  \mu^{h,\dl} (\xi^{h,\dl}(s), \al^{h,\dl}(s), \phi) m^{h,\dl}_s(d\phi) ds \\
&\quad\disp+ \int_{0}^{t}  \int_{\U} \sg^{h,\dl}  (\xi^{h,\dl}(s), \al^{h,\dl}(s), \phi ) m^{h,\dl}_s(d \phi) dW^{h,\dl}(s)\\
&\quad\disp +\bigg(\int_0^t \int_{\rr_+} \int_\U \wdt q^h_1 (q_1( \al^{h,\dl}(s),Z^{h,\dl}(s), \rho_1), a_1)  m^{h,\dl}_s(d \phi) N^{h,\dl}_1(ds, d\rho_1) \\
&\qquad\disp - \int_0^t \int_{\rr_+} \int_\U \ka_1 \wdt q^h_2(q_2(\al^{h,\dl}(s), Z^{h,\dl}(s), \rho_2), a_2)m^{h,\dl}_s(d\phi) N^{h,\dl}_2(ds,d\rho_2)\bigg) e_1\\
&\quad \disp+\bigg(\int_0^t \int_{\rr_+} \int_\U \wdt q^h_2 (q_2( \al^{h,\dl}(s),Z^{h,\dl}(s), \rho_2), a_2)  m^{h,\dl}_s(d \phi) N^{h,\dl}_2(ds, d\rho_2) \\
&\qquad\disp - \int_0^t \int_{\rr_+} \int_\U \ka_2 \wdt q^h_1(q_1(\al^{h,\dl}(s), Z^{h,\dl}(s), \rho_1), a_1)m^{h,\dl}_s(d\phi) N^{h,\dl}_1(ds,d\rho_1)\bigg) e_2 + \e^{h,\dl}(t).
\earray\eeq
We can also rewrite \eqref{mform} as
\beq{mform1}\barray
X(t) &\disp= x+ \int_0^t \int_\U \mu(X(s), \al(s), \phi(s)) m_s(d\phi) ds +\int_0^t \int_\U \sg(X(s), \al(s),\phi) m_s(d \phi) dW(s)\\
&\disp\quad + \bigg( \int_0^t \int_{\rr_+} \int_\U \wdt q_1 (q_1(\al(s), Z(s), \rho_1), a_1) m_s(d\phi) N_1(ds, d\rho_1)\\
&\disp\qquad - \int_0^t \int_{\rr_+} \int_\U \ka_1 \wdt q_2 (q_2(\al(s),Z(s), \rho_2), a_2) m_s(d\phi) N_2(ds, d\rho_2)\bigg) e_1\\
&\disp\quad +\bigg( \int_0^t \int_{\rr_+} \int_\U \wdt q_2 (q_2(\al(s),Z(s), \rho_2), a_2) m_s(d\phi) N_2(ds, d\rho_2) \\
&\disp\qquad - \int_0^t \int_{\rr_+} \int_\U \ka_2 \wdt q_1 (q_1(\al(s), Z(s), \rho_1), a_1) m_s(d\phi) N_1(ds, d\rho_1) \bigg) e_2.
\earray\eeq

\subsection{Convergence of Approximating Markov Chains}

\begin{lem} \label{mc}
Using the transition probability $\{p^{h, \dl}\cd \}$ defined in \eqref{prob2}, the interpolated process of the constructed Markov chain $\{ \al^{h,\dl} \cd$ converges weakly to $\al \cd$, the Markov chain with generator $Q$.
\end{lem}

\para{Proof.} The proof can be obtained similar to \cite[Theorem 3.1]{YinZB03}. The details are thus omitted. \qed

\begin{thm}\label{tight} Let the approximating chain $\{ \xi^{h, \dl}_n, \al^{h, \dl}_n, n < \infty\}$ constructed with transition probabilities defined in \eqref{prob2} be locally consistent with \eqref{cst}, $m^{h, \dl} \cd$ be the relaxed control representation of $\{u^{h, \dl}, n < \infty\}$, $(\xi^{h, \dl} \cd, \al^{h, \dl} \cd)$ be the continuous-time interpolation defined in \eqref{ipl}. Then $\{ \xi^{h, \dl} \cd, \al^{h, \dl} \cd, m^{h, \dl} \cd, W^{h, \dl} \cd, N^{h, \dl}_1 (\cdot, \cdot), N^{h, \dl}_2 (\cdot, \cdot) \}$ is tight.
\end{thm}

\para{Proof.}
Note that $\al^{h,\dl} \cd$ is tight.
It follows that
for each $\Dl >0$, each $t>0$, and $0< \wdt t \le \Delta$,
there is a random variable $\gamma^{h,\delta} (\Delta) >0$ such that
\beq{Var}\barray
\ex_t |W^{h,\dl} (t+ \wdt t) - W^{h,\dl}(t)|^2 &= \disp\sum_{\lf t/\dl \rf}^{ \lf t+\wdt t/\dl \rf-1} \ex_t \{[ \sg^{h,\dl} (x,\phi, i) ]^{-1}(\Dl\xi^{h,\dl}_r -\exr \Dl \xi^{h,\dl}_r) I_{H^{h,\dl}_r} \}^2 \\
\leq \ex_t \gamma^{h,\dl} (\Delta)
\earray\eeq
satisfying
$
\disp\lim_{\Dl \ar 0} \limsup_{ h,\dl \ar 0} \ex \gamma^{h,\dl}(\Delta) =0,
$
which
 yields the tightness of $W^{h,\dl} \cd$. A similar argument leads to  the tightness of $M^{h,\dl} \cd$. The sequence $m^{h,\dl} \cd$ is tight because
 of its compact range space. By virtue of
 \cite{KushnerD}[Theorem 9.2.1], we obtain the tightness of $\{ N_k^{h,\dl}\cd ,\ k =1,2 \}$ since the mean number of claims on any bounded interval $[t, t + t_1]$ is bounded and $$\disp \lim_{\Dl \ar 0} \inf_{r} \pr \{\nu^{h,\dl}_{k, r+1} - \nu^{h,\dl}_{k, r} > \Dl | \nu^{h,\dl}_{k,r}\}=1.$$ This implies the tightness of $\{R_k^{h,\dl} \cd, \ k =1,2 \}$. As a consequence, $\xi^{h,\dl} \cd$ is tight and $$\{ \xi^{h,\dl} \cd, \al^{h,\dl} \cd, m^{ h,\dl} \cd, W^{h, \dl} \cd, N^{ h,\dl}_1 (\cdot, \cdot), N^{h,\dl}_2 (\cdot, \cdot) \} \ \hbox{ is tight.}$$ \qed

Because $\{ \xi^{h,\dl} \cd, \al^{h,\dl} \cd, m^{ h,\dl} \cd, W^{h, \dl} \cd, N^{ h,\dl}_1 (\cdot, \cdot), N^{h,\dl}_2 (\cdot, \cdot) \}$ is tight, the Prohorov's theorem implies that it is sequentially compact. Thus we can extract a weakly convergent subsequence. Select such a convergent subsequence and still index the sequence by $h,\dl$ for notational simplicity. We proceed to characterize the limit process.

\begin{thm} Let $\{ \xi \cd, \al \cd, m \cd, W \cd, N_1 (\cdot, \cdot), N_2 (\cdot, \cdot) \}$ be the limit of weakly convergent subsequence and $\F_t$ be the $\sg$-algebra generated by $\{X(s), \al(s), m(s), W(s), N_1(s, \cdot), N_2(s, \cdot), s \le t \}$. Then $W\cd$ is a standard $\F_t$-Brownian motion and $N_1( \cdot, \cdot)$, $N_2(\cdot, \cdot)$ are $\F_t$-Poisson measures, and $m \cd$ is an admissible relaxed control.
\end{thm}

\para{Proof.}
The proof is divided into several steps.

\underline{Step 1}:
By the Skorohod representation, with a slight abuse of notation, $\{\xi^{h,\dl} \cd, \al^{h,\dl} \cd, m^{h,\dl} \cd$, $W^{h,\dl} \cd, N^{h,\dl}_1 (\cdot, \cdot), N^{h,\dl}_2 (\cdot, \cdot)\}$ converges to $\{\xi \cd, \al \cd, m \cd, W \cd, N_1(\cdot, \cdot), N_2 (\cdot, \cdot) \}$ w.p.1, and the convergence is uniform on any compact set.

To proceed,  we first verify that $W \cd$ is an $\F_t$-Brownian motion. For any real-valued and continuous function $\psi$, define
\beq{vc}
(\psi, m)_t = \int_0^t \int_\U \psi(\phi, s)m_s(d \phi).
\eeq
For any given $f\cd \in {\cal C}^2_0 (\rr^3)$
($C^2$ function with compact support), consider an associate operator
$\disp \Ll_w f(w) = \frac{1}{2} \sum_{i=1}^3 \dfrac{\partial^2}{\partial w_i \partial w_i} f(w)$.
Let $t, \wdt t >0$ be given with $t + \wdt t \leq T$, along with arbitrary positive integers $\ka$ and $\wdt \ka$, arbitrary $t_i \leq t$ and continuous functions $\psi_j$ with $i \leq \ka$ and $j \leq \wdt \ka$, any bounded and continuous function $h \cd$, and arbitrary $f \in {\cal C} ^2_0(\rr^3)$. Denote $\{\Ga^\ka_i, i \le \ka \}$ as a sequence of nondecreasing partition of $\rr_+$ such that $\Pi (\partial \Ga^\ka_i)=0$ for all $i, \ka$, where $\partial \Ga^\ka_i$ is the boundary of the set $\Ga^\ka_i$. As $\ka \ar \infty$, let the diameter of the sets $\Ga^\ka_i$ go to zero.

By \eqref{BM}, $W^{h,\dl} \cd$ is an $\F^{h,\dl}_t$-Brownian motion,
\bea
&\ex h(\xi^{h,\dl}(t_i),\al^{\dl}(t_i), W^{h,\dl}(t_i), (\psi_j, m^{h,\dl})_{t_i}: i \le \ka, j \le \wdt \ka) \\
&\times \displaystyle{ \bigg(f(W^{h,\dl}(t+ \wdt t) - f(W^{h,\dl}(t)) - \int_t^{t + \wdt t} \Ll_w f(W^{h,\dl} (s)) ds  \bigg) =0.}
\eea
By the weak convergence and the Skorohod representation,
 we may assume that $W^{h, \dl} \cd$ converges to $W \cd$ w.p.1, and hence as $\dl \ar 0$,
\bea
&\disp \ex h(\xi^{h,\dl}(t_i),\al^{\dl}(t_i), W^{h,\dl}(t_i), (\psi_j, m^{h,\dl})_{t_i}: i \le \ka, j \le \wdt \ka) \bigg(f(W^{h, \dl}(t+ \wdt t) - f(W^{h,\dl}(t))  \bigg)\\
&\disp\qquad \ar \ex h(\xi(t_i),\al(t_i), W(t_i), (\psi_j, m)_{t_i}: i \le \ka, j \le \wdt \ka) \bigg(f(W(t+ \wdt t) - f(W(t))  \bigg), \\
&\disp \ex h(\xi^{h,\dl}(t_i), \al^{h,\dl}(t_i), W^{h,\dl}(t_i), (\psi_j, m^{h,\dl})_{t_i}: i \le \ka, j \le \wdt \ka) \bigg(\int_t^{t + \wdt t} \Ll_w f(W^{h,\dl} (s)) ds \bigg) \\
& \disp\qquad \ar \ex h(\xi(t_i),\al(t_i), W(t_i), (\psi_j, m)_{t_i}: i \le \ka, j \le \wdt \ka) \bigg(  \int_t^{t + \wdt t} \Ll_w f(W (s)) ds \bigg).
\eea
Thus,
\bea
\disp\ex h(\xi(t_i), \al(t_i), W(t_i), (\psi_j, m)_{t_i}: i \leq \ka, j \le \wdt \ka ) \bigg(f(W(t + \wdt t)) - f(W(t))  - \int_t^{t + \wdt t} \Ll_w f(W(s)) ds\bigg) = 0.
\eea
Moreover, consider the collection of random variables generated by $\{ h(\xi(t_i),\al(t_i), W(t_i), (\psi_j, m)_{t_i}: i \leq \ka, j \le \wdt \ka) \}$.
It follows that $f(W(t))- f(W(0)- \int^t_0 \Ll_w f(W(s))ds$  is a martingale.
 By virtue of
 the Skorohod representation and the dominated convergence theorem together with \eqref{Var}, we have
\bea
\ex h(\xi(t_i), \al(t_i), W(t_i), (\psi_j,m)_{t_i}: i \leq \ka, j \le \wdt \ka) \bigg[(W(t+ \dl)- W(t))^2 - \dl \bigg] = 0.
\eea
Moreover,
the quadratic variation of the martingale $W(t)$ is $tI_3$ and $W \cd$ is an $\F_t$-Brownian motion, where $I_3$ is an $3 \times 3$ identity matrix.

\underline{Step 2}:
We proceed
to show that $N_k(\cdot, \cdot)$ is an $\F_t$-Poisson measure for each $k=1,2$. Let $\theta \cd$ be a continuous function on $\rr_+$ and define the process
\bea
\disp\Theta_k = \int_0^t \int_{\rr^+} \theta (\rho) N_k(ds,d\rho).
\eea
Using similar
argument as in the proof of
the Brownian motion above, if $f \cd \in {\cal C}_0^2 (\rr^3)$ then
\bea
&\ex h(\xi(t_i),\al(t_i), W(t_i), (\psi_j,m)_{t_i}, N(t_i, \Ga^\ka_i), i \le \ka, j \le \wdt \ka) \\
&\quad \times \bigg[ f(\Theta_k (t + \wdt t)) - f(\Theta_k(t)) - \la_k \int_t^{t + \wdt t} \displaystyle{\int_{\rr_+} \big[ f(\Theta_k(s) + \theta(\rho)) - f(\Theta_k(s)) \big] \Pi(d \rho) ds}\bigg] =0.
\eea
This implies that $N_k(\cdot, \cdot)$ is an $\F_t$-Possion measure for each $k=1,2$.

\underline{Step 3.1}: We will use \eqref{mc1} for the rest of the proof. Note that $\disp \ex |\e^{h,\dl}(t)| \ar 0$ as $h,\dl \ar 0$. Letting $h,\dl \ar 0$ and using the Skorohod representation for \eqref{mc1}, we have
\bea
\disp\int_0^t \int_\U \mu^{h, \dl}(\xi^{h,\dl}(s), \al^{h,\dl}(s), \phi) m^{h,\dl}_s(d \phi) ds - \int_0^t \int_\U \mu(\xi(s), \al(s), \phi) m^{h,\dl}_s(d\phi) ds \ar 0,
\eea
uniformly on any bounded time interval with probability one. On the other hand, the sequence $m^{h,\dl} \cd$ converges in the compact-weak topology, thus, for any continuous and bounded function $ \psi \cd$ with compact support,
\bea
\disp\int_0^t \int_\U \psi(\phi,s) m^{h,\dl} (d\phi ds) \ar \int_0^t \int_\U \psi(\phi, s) m(d\phi ds) \hbox{ as } h,\dl \ar 0.
\eea
By virtue of
the Skorohod representation and the weak convergence,
 as $h,\dl \ar 0$,
\beq{mcv}
\int_0^t \int_\U \mu^{h,\dl}(\xi^{h,\dl}(s), \al^{h,\dl}(s), \phi) m^{h,\dl}_s (d\phi) ds - \int_0^t \int_\U \mu(\xi(s),\al(s),\phi) m_s(d\phi) ds \ar 0,
\eeq
uniformly in $t$ with probability one on any bounded interval.

\underline{Step 3.2}: For any $t^1 \ge t$, $t^2 \ge 0$ with $t^1 + t^2 \le T$, any ${\cal C}_0^{1,2}$ function $f \cd$ (functions that have compact support whose first partial derivative w.r.t. the time variable and the second partial derivatives w.r.t. the state variable $x$ are continuous), bounded and continuous function $h \cd$, any positive integers $\ka, \wdt \ka,$ $t_i$, and any continuous function $\psi_j$  satisfying $t \le t_i \le t_1$ and $i \le \ka$, and $j \le \wdt \ka$, the weak convergence and the Skorohod representation imply that
\beq{cv1}\barray
&\ex h(\xi^{h,\dl}(t_i), \al^{h,\dl}(t_i),M^{h, \dl}(t_i),(\psi_j,m^{h, \dl})_{t_i}: i \le \ka, j \le \wdt \ka)
\\ &\quad \times [f(t^1+t^2, M^{h, \dl}(t^1 +t^2)) - f(t^1, M^{h, \dl}(t^1))] \\
&\ar \ex h(\xi(t_i), \al(t_i),M(t_i),(\psi_j, m)_{t_i}: i \le \ka, j \le \wdt \ka) [f(t^1+t^2, M(t^1+t^2) - f(t^1, M(t^1))]
\earray\eeq
with  $h,\dl \ar 0.$
Choose a sequence $\{ n^\dl\}$ such that $n^\dl \ar \infty$ but $\Dl^\dl = \dl n^\dl \ar 0$, then
\beq{sum}\barray
& \ex h(\xi^{h,\dl}(t_i), \al^{h,\dl}(t_i),M^{h, \dl}(t_i), (\psi_j,m^{h,\dl})_{t_i}: i \le \ka, j \le \wdt \ka)\\
& \quad \times [f(t^1+t^2, M^{h, \dl}(t^1+t^2)) - f(t^1, M^{h, \dl}(t^1))] \\
& =\disp\ex h(\xi^{h,\dl}(t_i), \al^{\dl}(t_i),M^{h,\dl}(t_i),(\psi_j, m^{h,\dl})_{t_i}: i \le \ka, j \le \wdt \ka) \\
& \quad \times \bigg[\sum_{l n^\dl = t^1/\dl}^{(t^1+t^2)/\dl-1} f(\dl(l n^\dl + n^\dl), M^{h,\dl}(\dl(l n^\dl + n^\dl)))\\
&\disp\quad  - f(\dl ln^\dl, M^{h,\dl}(\dl( l n^\dl+n^\dl)))
+f(\dl l n^\dl, M^{h,\dl}(\dl (l n^\dl +n^\dl))) - f(\dl l n^\dl, M^{\dl}(\dl l n^\dl))  \bigg].
\earray\eeq
Note that
\bea
 & \disp\sum_{l n^\dl = t^1/\dl}^{(t^1+t^2)/\dl} [f(\dl(l n^\dl + n^\dl), M^{\dl}(\dl(l n^\dl + n^\dl))) - f(\dl ln^\dl, M^{h,\dl}(\dl( l n^\dl+n^\dl)))]\\
 &\disp\quad = \sum_{l n^\dl = t^1/\dl}^{(t^1+t^2)/\dl} \sum_{k= l n^\dl}^{l n^\dl + n^\dl -1} [f(\dl(k+1), M^{h, \dl}(\dl(l n^\dl + n^\dl))) - f(\dl k, M^{h,\dl}(\dl( l n^\dl+n^\dl)))]\\
 &\disp\quad = \sum_{l n^\dl = t^1/\dl}^{(t^1+t^2)/\dl} \dfrac{\partial f(\dl l n^\dl, M^{\dl}(\dl(l n^\dl + n^\dl)))}{\partial s} \Dl^{\dl} + o(1),
\eea
where $o(1) \ar 0$ in mean uniformly in $t$ as $\dl \ar 0$. Letting $\dl l n^\dl \ar s$ as $\dl \ar 0$, then $\dl (l n^\dl + n^\dl) \ar s$ since $\Dl^\dl = \dl n^\dl \ar 0$ as $\dl \ar 0$.
Then, by the weak convergence and the Skorohod representation, the continuity of $h \cd$, and the smoothness of $f \cd$ imply that
\beq{cvt} \barray
&\disp \ex h(\xi^{h,\dl}(t_i), \al^{h,\dl}(t_i),M^{h,\dl} (t_i),(\psi_j,m^{h,\dl})_{t_i}: i \le \ka, j \le \wdt \ka)\\
&\times  \sum_{l n^\dl}^{(t^1+t^2)/\dl} \bigg[f(\dl(l n^\dl + n^\dl), M^{h, \dl}(\dl(l n^\dl + n^\dl))) - f(\dl ln^\dl, M^{h,\dl}(\dl( l n^\dl+n^\dl))) \bigg]
\\ & \ar \ex h(\xi(t_i), \al(t_i),M(t_i),(\psi_j, m)_{t_i}: i \le \ka, j \le \wdt\ka) \bigg[\displaystyle{\int_{t^1}^{t^1+t^2} \dfrac{\partial f(s, M(s))}{\partial s} ds} \bigg]
\disp \qquad  \hbox{ as } h,\dl \ar 0.
\earray \eeq
The 
last part of \eqref{sum} can be seen as
\bea
&\disp \sum_{l n^\dl = t^1/\dl}^{(t^1+t^2)/\dl} [f(\dl l n^\dl, M^{h,\dl}(\dl(l n^\dl + n^\dl))) - f(\dl l n^\dl, M^{h,\dl}(\dl l n^\dl))] \\
&\disp  = \sum_{l n^\dl = t^1/\dl}^{(t^1+t^2)/\dl}  \bigg\{ \sum_{i=1}^3\frac{1}{2} f_{M_i M_i}(\dl l n^\dl, M^{h, \dl}(\dl l n^\dl)) \sum_{k = l n^\dl}^{l n^\dl + n^\dl -1} [M_i^{h,\dl}(\dl(l n^\dl +n^\dl)) - M^{h,\dl}_i (\dl l n^\dl)]^2 \bigg\}\\
&\disp \qquad  + \wdt \e^{h,\dl} (t^1+t^2) -\wdt \e^{h,\dl}(t^1),
\eea
where $f_{M_i M_i} $ denotes the second partial derivatives, $M^{h,\dl}_i \cd$ is the $i$-th component of $M^{h,\dl} \cd$ and $\disp\sup_{t \le t^1 \le T} \ex |\wdt \e^{h,\dl} (t^1)| \ar 0$ as $h,\dl \ar 0$.

By \eqref{eq:not} and the definition of $\sg^{h,\dl} \cd $, we have
\beq{cvv}\barray
 &\disp \sum_{l n^\dl = t^1/\dl}^{(t^1+t^2)/\dl}  \bigg\{\sum_{i=1}^2 f_{M_i M_i}(\dl l n^\dl, M^{h, \dl}(\dl l n^\dl)) \sum_{k = l n^\dl}^{l n^\dl + n^\dl -1} [M^{h,\dl}_i(\dl(l n^\dl +n^\dl)) - M^{h,\dl}_i (\dl l n^\dl)]^2 \bigg\} \\
 &\disp \ar \int_{t^1}^{t^1+t^2} \int_\U Tr\bigg[H_M f(s, M(s)) \sg(\xi(s), \al(s), \phi) [\sg(\xi(s), \al(s), \phi)]'\bigg] m_s(d \phi) ds,
\earray \eeq
where $H_M(f(s, M(s)))$ is the Hessian matrix of $f\cd$
at time $s$, $Tr \cd$ represents for the trace of a matrix.

Using \eqref{cv1}-\eqref{cvv}, we have
\bea
&\disp\ex h(\xi(t_i), \al(t_i), M(t_i),(\psi_j, m)_{t_i}: i \le \ka, j \le \wdt \ka)\\
&\disp \quad \times \bigg[f(t^1+t^2, M(t^1+t^2)) - f(t^1, M(t^1)) -\int_{t^1}^{t^1+t^2} \dfrac{\partial f(s,M(s))}{\partial s}ds \\
&\disp\qquad - \int_{t^1}^{t^1+t^2} \int_\U \frac{1}{2}Tr\bigg(H_M f(s, M(s)) \sg(\xi(s), \al(s), \phi) [ \sg(\xi(s), \al(s), \phi)]'\bigg)  m_s(d \phi) ds \bigg] =0.
\eea
Therefore,
\beq{mac}
\int_0^t \int_\U \sg^{h,\dl}(\xi^{h,\dl}(s), \al^{h, \dl} (s), \phi) m^{h,\dl}_s (d\phi) ds \ar \int_0^t \int_\U \sg(\xi(s), \al(s), \phi) m_s(d \phi) ds,
\eeq
uniformly in $t$ with probability one on any bounded interval.

Using the same arguments
as in Step 2 and Step 3.2, we obtain the
limits for the latter parts of \eqref{mc2}. As a result, $\xi\cd$ is the solution of \eqref{mform1}, which means $\xi(s) =X(s)$ w.p.1 and $m \cd$ is an admissible relaxed control. \qed

\subsection{Convergence of the Cost and the Value Functions}


Note that since $U_k \cd$ satisfies the Inada's conditions. There exist positive real numbers $K$ and $k_0$ such that $|U_k(X_k)| \le K(1 + |X_k|^{k_0})$. We proceed to prove the following result.

\begin{thm}
Suppose that the utility functions $U_k \cd$ has at most polynomial growth. Then the value functions $V^{k, h,\dl}(t, x, i)$
converges to $V^k(t,x,i)$ for $k=1,2$,
respectively, as $h,\dl \ar 0$.
\end{thm}

\para{Proof.}
 By
 \thmref{tight}, each sequence $\{\xi^{h,\dl}\cd, \al^{h,\dl} \cd, m^{h,\dl} \cd, W^{h,\dl} \cd, N^{h,\dl}_1 (\cdot, \cdot), N^{h,\dl}_2 (\cdot, \cdot) \}$ has a weakly convergent subsequence with the limit $\xi^{h,\dl} \cd$ satisfying \eqref{mform1}. Using the same notation as above and applying the Skorohod representation, the weak convergence,
 as $h,\dl \ar 0$,
 $J^{k,h,\dl} (t, x, i, m^{h,\dl}) \ar J^{k}(t,x,i,m )$,
 for $k=1,2$.
 The cost function is given by \eqref{dp}. Since $V^k (t, x, i)$ is the maximizing expected utility, for any admissible control $m \cd$,
 $J^{k}(t, x,i, m) \le V^k (t,x, i)$,
 for $k =1,2$.
 Let $\wdt m ^{h,\dl} \cd$ be an optimal relaxed control for $\{\xi^{h,\dl} \cd \}$, which
 implies
 \bea
 \disp V^{k,h,\dl} (t, x, i) = J^{k,h, \dl} (t, x, i, \wdt m^{h,\dl}) = \sup_{m^{h,\dl}} J^{k,h,\dl}(t,x, i,m^{h,\dl}).
 \eea
 Choose a subsequence $\{ \wdt h, \wdt \dl \}$ of $\{h,\dl \}$ such that
 \bea
 \disp  \limsup_{ h,\dl \ar 0} V^{k, h, \dl} (t, x,i) =\lim_{\wdt h, \wdt \dl \ar 0} V^{k,\wdt h,\wdt \dl} (t,x,i)= \lim_{\wdt h,\wdt \dl \ar 0} J^{k,\wdt h, \wdt \dl} (t,x,i, \wdt m^{\wdt h,\wdt \dl}).
 \eea
 Without loss of generality, we may assume that $\{\xi^{\wdt h, \wdt \dl}\cd, \al^{\wdt h, \wdt \dl} \cd, W^{\wdt h, \wdt \dl} \cd, \wdt m ^{\wdt h,\wdt \dl} \cd, N_1^{\wdt h,\wdt \dl} (\cdot, \cdot), N_2^{\wdt h, \wdt \dl} (\cdot, \cdot) \}$ converges weakly to $\{ X \cd, \al \cd, W \cd, m \cd, N_1(\cdot, \cdot), N_2(\cdot, \cdot)\}$, where $m \cd$ is an admissible relaxed control. Then the weak convergence and the Skorohod representation
 leads to
\bea
\disp\limsup_{ h,\dl \ar 0} V^{k, \dl}(t,x,i) =\lim_{\wdt h,\wdt \dl \ar 0} J^{k, \wdt h, \wdt \dl} (t,x,i, \wdt m^{\wdt h,\wdt \dl})= J^k (t, x,i , m) \le V^k (t, x, i).
\eea
We claim that $\disp\liminf_{\dl} V^{k, \dl} (t,x,i) \ge V^k(t,x,i)$.

Suppose that $\lbar m \cd$ is an optimal control with Brownian motion $W \cd$ such that $\lbar X \cd$ is the associated trajectory. By 
the chattering lemma (see \cite{YinJJ09} and page 59-60 of \cite{KushnerD}),
for any given $\eta, \dl_\eta >0$, there is an $\e >0$ and an ordinary control $u^{\eta, \dl_\eta} \cd$ that takes only finite many values, $u^{\eta, \dl_\eta} \cd$ is a constant in $[ \iota \e, \iota \e +\e)$, $\lbar m^{\eta, \dl_\eta} \cd$ is its relaxed control representation, and $J^k(t,x,i, \lbar m^{\eta, \dl_\eta}) \ge V^k(t,x,i) -\eta$.
For each $\eta , \dl_\eta>0$, and the corresponding $\e >0$, consider an optimal control problem with piecewise constant on $[\iota \e, \iota \e + \e)$. We consider the process $\{X^{\eta, \dl_\eta} (\iota \e),\al^{\eta, \dl_\eta}, m^{\eta, \dl_\eta} (\iota \e), W^{\eta, \dl_\eta} (\iota \e) \}$.
Let $\hat{u}^{\eta,\dl_\eta} \cd$ be the optimal control, $\hat{m}^{\eta,\dl_\eta} \cd$ the relaxed control representation, and $\hat{X}^{\eta, \dl_\eta} \cd$ the associated trajectory. Since $\hat{m}^{\eta, \dl_\eta} \cd$ is the optimal control, $J^k(t,x,i, \hat{m}^{\eta, \dl_\eta}) \ge J^k(t,x,i, \lbar m ^{\eta, \dl_\eta}) \ge V^k(t,x,i) - \eta$.
Using the chattering lemma,
we can approximate $\hat{m}^{\eta, \dl_\eta} \cd$ by a sequence of $m^{h,\dl} \cd$. Then \bea V^{k,h,\dl} (t,x,i) \ge J^{k,h,\dl}(t,x,i, m^{h,\dl}) \ar J^k (t,x,i, \hat{m}^{\eta,\dl_\eta}).   \eea Moreover, \bea \disp \liminf_{h,\dl \ar 0} V^{k,h,\dl} (t,x,i) \ge \lim_{h,\dl \ar 0}J^{k,h,\dl}(t,x,i,m^{h,\dl}) = J^k(t,x,i,\hat{m}^{\eta,\dl_\eta}).\eea Thus, $\disp \liminf_{h,\dl \ar 0} V^{k,h,\dl} (t,x,i) \ge V^k(t,x,i) - \eta$. The arbitrariness of $\eta$ implies that $\disp \liminf_{\dl \ar 0} V^{k,\dl} (t,x,i) \ge V^k(t,x,i)$, which completes the proof.
\qed

\section{Numerical Examples}\label{sec:exam}
In this section, we present some numerical results for the case in which both insurance companies are constant absolute risk aversion (CARA) agents, i.e., each agent has an exponential utility function. More precisely, the utility function of each insurer has the form
\beq{ut}
U_k(X_k) = -\frac{1}{\eta_k} \exp{(-\eta_k X_k)}, \quad \hbox{for } \eta_k>0,\ k =1,2.
\eeq
Based on the algorithm constructed above, we carry out the computation by valuing iterations in a backward manner by time.
\begin{itemize}
\item[1.] Set $t =T-\dl$ and $J^{k, h, \dl}(T, x, i, u^{h,\dl}) = U_k(x_k)$, for each $k =1,2$.

\item[2.] By \eqref{dp}, we  obtain
\bea
J^{k, \dl}(t, x, i, u^{h,\dl}) &= \disp(1- \la \dl ) \sum_{(y,j)} p^{h,\dl}_k\bigg((x,i),(y,j)|u^{h,\dl} \bigg) J^{k,\dl}(t+\dl, y,j,u^{h,\dl})\\
&+ \disp \la_k \dl \int_{\rr_+}J^{k,\dl}(t+\dl, x_k -\wdt q_{k}^h(q_k(i,z,\rho_k), a^{h,\dl}_k), z, i, u^{h, \dl}) \Pi_k(d \rho_k)\\
&+ \disp \la_l \dl  \int_{\rr_+}J^{k,\dl} (t+ \dl, x_k + \kappa_k \wdt q_l^h(q_l(i,z,\rho_l), a^{h,\dl}_l),z,i,u^{h,\dl}) \Pi_l(d \rho_l).
\eea
Find the pair $\{ \hat{u}^{h,\dl}_k , k=1,2 \} $  and record $u^{h,\dl}_k (t) = \hat{u}^{h,\dl}_k $ satisfying that for any $u^{h,\dl}_k \in \U_k$,:
\bea
& J^{1,h, \dl} (t, x, i, u^{h,\dl}_1,  \hat{u}^{h,\dl}_2) \le J^{1,h, \dl} (t, x, i, \hat{u}^{h,\dl}_1,  \hat{u}^{h,\dl}_2) \\
& J^{2,h, \dl} (t, x, i, \hat{u}^{h,\dl}_1,  u^{h,\dl}_2) \le J^{2,h, \dl} (t, x, i, \hat{u}^{h,\dl}_1,  \hat{u}^{h,\dl}_2).
\eea

\item[3.] Let $t = t-\dl$ and continue the procedure until $t =0$.
We consider the case in which the discrete event consists of two states, or equivalently, the Markov chain has two states with given claim size distributions. In addition, we assume that the claim size distributions are identical in each regime. By using the value iteration methods, we numerically solve the optimal control problems.
The continuous-time Markov chain $\al(t)$ representing the discrete event state has the generator $Q = \begin{pmatrix}
&-0.5 & 0.5\\
&0.5&-0.5
\end{pmatrix}$ and takes values in $\M = \{ 1,2\}$.
\end{itemize}

The parameters of the utility function are $\eta_1=17.0$ and $\eta_1=21.0$ respectively. The sensitivities are $ \kappa_1= 0.8$ and $ \kappa_2= 0.7$. The claim severity of both players follows exponential distribution $f\left(y\right) = \theta_k e^{-\theta_k y}$ with $\theta_1=0.3$ and $\theta_2=0.2$. To incorporate the difference between claim densities in different regimes, we assume arriving rates of Poisson jump are different. So is the setup for risk-free return, and premium income rate. The detail of setup is as follows in Table \ref{table:par}.
	
\begin{figure}[htbp!]
		\centering
		\includegraphics[scale=0.49]{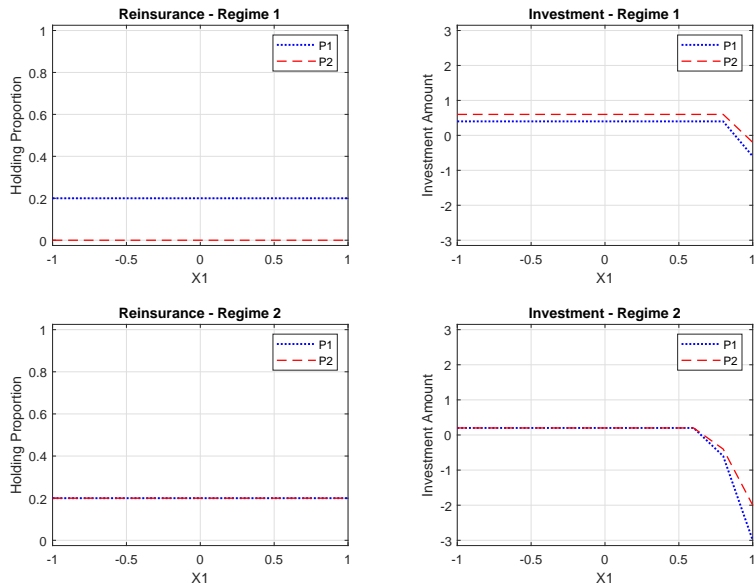}
		\caption{Controls for varying $X_1$ with $T=0.08,$ $Z=1.01,$ and $X_2=0$.}
		\label{P1}
	\end{figure}

\begin{figure}[htbp!]
	\centering
	\includegraphics[scale=0.49]{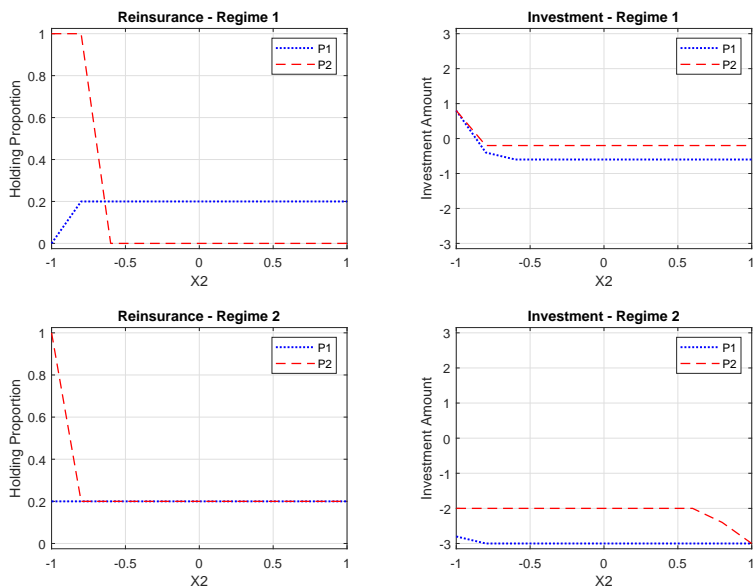}
	\caption{Controls for varying $X_2$ with $T=0.08,$ $Z=1.01,$ and $X_1=1$.}
	\label{P2}
\end{figure}

	\begin{table}[htbp]
		\centering
		\begin{tabular}{c c c c c  c c  }
			\hline
			Regime & $r$ & $c_1$  &$c_2$ & $\lambda_1$  &$\lambda_2$ &\\
			\hline
			&0.02&0.05&0.02&0.20&0.30\\
			2&0.03&0.10&0.20&0.80&0.70\\
			\hline
		\end{tabular}
		\caption{Parameters values}\label{table:par}
	\end{table}
The reinsurance premium rates are computed from expectation premium principle as:
\begin{equation}\nonumber
g_k\left(a_k\right)=\left(1+l_k\right)\left(1-a_k\right)\mathbb{E}\left[A_k\right],
\end{equation}
where $l_1=1.1$ and $l_2=1.15$. Further, the volatility and the drift of the financial market index and the risky return rate are modeled respectively by:
\begin{equation}\nonumber
\begin{split}
\mu_S\left(i,Z\right)&=0.2iZ;\\
\sigma_S\left(i,Z\right)&=0.4iZ;\\
\mu_Z\left(t,Z\right)&=0.4\left(t+1\right)Z;\\
\sigma_Z\left(t,Z\right)&=0.1\left(t+1\right)Z.\\
\end{split}
\end{equation}
State discretization follows $\delta = 0.04$ and $h=0.2$. For reinsurance, we discretize the reinsurance rate into six levels uniformly located from 0 to 1. The investment amount is free from restrictions, and can vary from -3 to 3 with 0.2 increase.  The trend of investment and reinsurance for the varying relative surplus of player one is plotted in Figure \ref{P1}, and that for the second player is in Figure \ref{P2}.

From Figure \ref{P1}, we can observe that both players always hold a low proportion of claim, which is due to the fact that both players are very risk-averse. A big claim will not only reduce their relative surplus but also drive the surplus of their opponent side up. Precisely, in regime one, the proportion held by player two is 0, which is less than 0.2 of player one. This results from the claim arriving rate of the player two is relatively much higher than that of player one. In regime two, considering the high premium income rate and the same expected claim amount, a small proportion of claim is affordable for both players.

In view of the investment part of Figure \ref{P1}, player two's investment amount is always no less than that of  player one. Since the relative sensitivity of  player one to player two is higher, and player two's relative surplus is at a higher relative level initially, player two is willing to accept more risk for higher expected return in order to beat player one. Then, along with growing $X_1$, since player one's condition is improved, he tends to bear less risk. Meanwhile, player two adopts conservative strategy as well, since he can lower the uncertainty and make use of his advantage that he owns a higher premium income rate. Because a higher market volatility in regime two will introduce more risk, both players choose the investment amount much closer to 0 to lower the uncertainty.

Similar results can be seen
from Figure \ref{P2}.  For the reinsurance part, player two initially holds full proportion of a claim, since $X_2=-1$ leaves player two in a relative bad situation compared to $X_1=1$. To change this situation, he chooses to bear lots of risk to reduce the loss from reinsurance premium. This is more obvious for regime 1, where the claim arriving rate is lower.  Hence, in more risky scenarios where claim arriving rate is higher, players show relative risk averse and transfer more risks by reinsurance tools.
For the investment part, we can see that player two hold higher positions in risky assets in his portfolio in both regimes, which is consistent with the observations in Figure \ref{P1}.

\section{Concluding Remarks}\label{sec:rem}

In this paper, we considered a non-zero-sum stochastic investment and reinsurance game between two insurance companies.  Both proportional and non-proportional reinsurance contracts were considered. Although we are able to obtain the systems of HJI equations using dynamic programming principle, solving the problem explicitly is
virtually
impossible. Thus we developed a numerical scheme using the Markov chain approximation method (MCAM) to solve the problem. Due to the complexity of the stochastic game formulation, even numerically  solving the systems of HJI equations is much more difficult than that of the previous work in stochastic optimization problems. The difficulties arise from the following two aspects. (1) With complex nonlinear state processes, the formulated high-dimension problem adds much difficulties in building approximating Markov chain. (2) The curse of dimensionality makes a significant impact
 and slow down the computation
 due to the large numbers of control variables and the dimensions of the HJI systems.
 Although the paper was devoted to a problem arising in risk management and insurance fields, the game problem formulation and the numerical methods developed can be more widely used in various other control and game problems.

For our problem,
the nature of the Markov chain approximation relies on building a high dimensional lattice of both driving state and control strategy to approximate the
value functions under different control scenarios. The optimization on every state follows the same computing rule, leading to the possibility
of using parallel acceleration techniques.
The first option coming to our mind is to incorporate multi-thread programming techniques into our completed C++ MCAM template library, which enables us to reduce development time by reusing the algorithm architecture of single-thread library. The latest eighth generation Intel CPUs are equipped with six computation cores, which allow maximal twelve threads to run simultaneously. If we parallelize the algorithm using ten threads, we can
enhance the time efficiency ten times. However, this is not enough to handle the computational complexity required for our problem. The high dimensionality requires the lattice to be very precise,
thus obtaining accurate results relies on generating a large number of nodes. Ten times acceleration seems a big enhancement, but it can only allow us to explore $10^\frac{1}{7}\approx 1.39$ times of $Z$, $X_1$, $X_2$, $a_1$, $a_2$, $b_1$, and $b_2$, which is
unable to meet requirements of the computational complexity.

GPU acceleration, e.g., CUDA, is another attractive choice here.
Although the frequency of GPU core is much lower than that of CPU, the number of GPU cores is usually hundreds of times of the number of CPU cores, and this makes GPU more suitable for parallel computation. The tenth generation NVidia GPU owns more than two thousand CUDA cores, which make it an easy solution for solving the complexity issue of MCAM algorithm on common stochastic optimization problem, where the maximal or the minimal value on a state is acquired from repeatedly comparing the newly computed value function value against the temporary optimal value so far. However, focusing on MCAM algorithm on our high dimensional game problem, CUDA acceleration is of very limited use. Not like CPU memory, which can be easily more than 64GB, the capacity of GPU memory is usually less than 8GB. The equilibrium strategy is obtained by searching on the value function information stored for different values of the control strategy. As a result, this memory consumption will occur for every GPU thread, which will easily lead the aggregated memory consumed by MCAM algorithm to exceed the GPU’s memory capacity.
From the above considerations, it appears that using
parallel programming techniques to high dimensional game problems 
needs a lot of
more thinking and effort.
Finding more efficient way for the numerical solution is our on-going work.

\end{document}